
\magnification=\magstep1
\hfuzz=1pt
\documentstyle{amsppt}
\let\fr\frak
\predefine\Sec{\S}
\define\C{{\Bbb C}}
\define\R{{\Bbb R}}


\define\Z{{\Bbb Z}}
\define\A{\text{\bf A}}

\define\G{\text{\bf G}}

\define\M{\text{\bf M}}
\define\N{\text{\bf N}}
\define\PPP{\text{\bf P}}


\define\z{\fr z}

\redefine\b{\beta}
\redefine\d{\delta}
\define\e{\varepsilon}
\define\ga{\gamma}
\predefine\l{\ll}
\define\k{\kappa}
\redefine\l{\lambda}
\redefine\o{\omega}

\define\s{\sigma}

\define\DE{\Delta}

\redefine\S{\Sigma}
\redefine\L{\Lambda}
\redefine\O{\Omega}
\define\T{\Theta}


\mathsurround=3pt

\define\diag{\operatorname{diag}}

\define\bs{\backslash}

\define\wt{\widetilde}

\define\Ind{\operatorname{Ind}}
\define\Res{\operatornamewithlimits{Res}}
\define\Hom{\operatorname{Hom}}

\message {REFERENCE FORMATTING FILES:}

\def\Resetstrings{
    \def\present{ }\let\bgroup={\let\egroup=}
    \def\Astr{}\def\astr{}\def\Atest{}\def\atest{}%
    \def\Bstr{}\def\bstr{}\def\Btest{}\def\btest{}%
    \def\Cstr{}\def\cstr{}\def\Ctest{}\def\ctest{}%
    \def\Dstr{}\def\dstr{}\def\Dtest{}\def\dtest{}%
    \def\Estr{}\def\estr{}\def\Etest{}\def\etest{}%
    \def\Fstr{}\def\fstr{}\def\Ftest{}\def\ftest{}%
    \def\Gstr{}\def\gstr{}\def\Gtest{}\def\gtest{}%
    \def\Hstr{}\def\hstr{}\def\Htest{}\def\htest{}%
    \def\Istr{}\def\istr{}\def\Itest{}\def\itest{}%
    \def\Jstr{}\def\jstr{}\def\Jtest{}\def\jtest{}%
    \def\Kstr{}\def\kstr{}\def\Ktest{}\def\ktest{}%
    \def\Lstr{}\def\lstr{}\def\Ltest{}\def\ltest{}%
    \def\Mstr{}\def\mstr{}\def\Mtest{}\def\mtest{}%
    \def\Nstr{}\def\nstr{}\def\Ntest{}\def\ntest{}%
    \def\Ostr{}\def\ostr{}\def\Otest{}\def\otest{}%
    \def\Pstr{}\def\pstr{}\def\Ptest{}\def\ptest{}%
    \def\Qstr{}\def\qstr{}\def\Qtest{}\def\qtest{}%
    \def\Rstr{}\def\rstr{}\def\Rtest{}\def\rtest{}%
    \def\Sstr{}\def\sstr{}\def\Stest{}\def\stest{}%
    \def\Tstr{}\def\tstr{}\def\Ttest{}\def\ttest{}%
    \def\Ustr{}\def\ustr{}\def\Utest{}\def\utest{}%
    \def\Vstr{}\def\vstr{}\def\Vtest{}\def\vtest{}%
    \def\Wstr{}\def\wstr{}\def\Wtest{}\def\wtest{}%
    \def\Xstr{}\def\xstr{}\def\Xtest{}\def\xtest{}%
    \def\Ystr{}\def\ystr{}\def\Ytest{}\def\ytest{}%
}
\Resetstrings

\def\Refformat{
         \if\Jtest\present
             {\if\Vtest\present\journalarticleformat
                  \else\conferencereportformat\fi}
            \else\if\Btest\present\bookarticleformat
               \else\if\Rtest\present\technicalreportformat
                  \else\if\Itest\present\bookformat
                     \else\otherformat\fi\fi\fi\fi}

\def\Rpunct{
   \def\Lspace{ }%
   \def\Lperiod{ }
   \def\Lcomma{ }
   \def\Lquest{ }
   \def\Lcolon{ }
   \def\Lscolon{ }
   \def\Lbang{ }
   \def\Lquote{ }
   \def\Lqquote{ }
   \def\Lrquote{ }
   \def\Rspace{}%
   \def\Rperiod{.}
   \def\Rcomma{,}
   \def\Rquest{?}
   \def\Rcolon{:}
   \def\Rscolon{;}
   \def\Rbang{!}
   \def\Rquote{'}
   \def\Rqquote{"}
   \def\Rrquote{`}
   }

\def\Lpunct{
   \def\Lspace{}%
   \def\Lperiod{\unskip.}
   \def\Lcomma{\unskip,}
   \def\Lquest{\unskip?}
   \def\Lcolon{\unskip:}
   \def\Lscolon{\unskip;}
   \def\Lbang{\unskip!}
   \def\Lquote{\unskip'}
   \def\Lqquote{\unskip"}
   \def\Lrquote{\unskip`}
   \def\Rspace{\spacefactor=1000}%
   \def\Rperiod{\spacefactor=3000}
   \def\Rcomma{\spacefactor=1250}
   \def\Rquest{\spacefactor=3000}
   \def\Rcolon{\spacefactor=2000}
   \def\Rscolon{\spacefactor=1250}
   \def\Rbang{\spacefactor=3000}
   \def\Rquote{\spacefactor=1000}
   \def\Rqquote{\spacefactor=1000}
   \def\Rrquote{\spacefactor=1000}
   }

\def\Refstd{
     \def\Acomma{\unskip, }
     \def\Aand{\unskip\ and }
     \def\Aandd{\unskip\ and }
     \def\Ecomma{\unskip, }
     \def\Eand{\unskip\ and }
     \def\Eandd{\unskip\ and }
     \def\acomma{\unskip, }
     \def\aand{\unskip\ and }
     \def\aandd{\unskip\ and }
     \def\ecomma{\unskip, }
     \def\eand{\unskip\ and }
     \def\eandd{\unskip\ and }
     \def\Namecomma{\unskip, }
     \def\Nameand{\unskip\ and }
     \def\Nameandd{\unskip\ and }
     \def\Revcomma{\unskip, }
     \def\Initper{.\ }
     \def\Initgap{\dimen0=\spaceskip\divide\dimen0 by 2\hskip-\dimen0}%
   }

\def\Smallcapsaand{
     \def\Aand{\unskip\bgroup{\Smallcapsfont\ AND }\egroup}%
     \def\Aandd{\unskip\bgroup{\Smallcapsfont\ AND }\egroup}%
     \def\eand{\unskip\bgroup\Smallcapsfont\ AND \egroup}%
     \def\eandd{\unskip\bgroup\Smallcapsfont\ AND \egroup}%
   }

\def\Smallcapseand{
     \def\Eand{\unskip\bgroup\Smallcapsfont\ AND \egroup}%
     \def\Eandd{\unskip\bgroup\Smallcapsfont\ AND \egroup}%
     \def\aand{\unskip\bgroup\Smallcapsfont\ AND \egroup}%
     \def\aandd{\unskip\bgroup\Smallcapsfont\ AND \egroup}%
   }

   \def\Citefont{}
   \def\ACitefont{}
   \def\Authfont{}
   \def\Titlefont{}
   \def\Tomefont{\sl}
   \def\Volfont{}
   \def\Flagfont{}
   \def\Reffont{\rm}
   \def\Smallcapsfont{\sevenrm}
   \def\Flagstyle#1{\hangindent\parindent\indent\hbox to0pt
       {\hss[{\Flagfont#1}]\kern.5em}\ignorespaces}


\def\Citebrackets{\Rpunct
   \def\Lcitemark{\def\Cfont{\Citefont}[\bgroup\Cfont}
   \def\Rcitemark{\egroup]}
   \def\LAcitemark{\def\Cfont{\ACitefont}\bgroup\ACitefont}%
   \def\RAcitemark{\egroup}
   \def\LIcitemark{\egroup}
   \def\RIcitemark{\bgroup\Cfont}
   \def\Citehyphen{\egroup--\bgroup\Cfont}
   \def\Citecomma{\egroup,\hskip0pt\bgroup\Cfont}%
   \def\Citebreak{}
   }

\def\Citeparen{\Rpunct
   \def\Lcitemark{\def\Cfont{\Citefont}(\bgroup\Cfont}
   \def\Rcitemark{\egroup)}
   \def\LAcitemark{\def\Cfont{\ACitefont}\bgroup\ACitefont}%
   \def\RAcitemark{\egroup}
   \def\LIcitemark{\egroup}
   \def\RIcitemark{\bgroup\Cfont}
   \def\Citehyphen{\egroup--\bgroup\Cfont}
   \def\Citecomma{\egroup,\hskip0pt\bgroup\Cfont}%
   \def\Citebreak{}
   }

\def\Citesuper{\Lpunct
   \def\Lcitemark{\def\Cfont{\Citefont}\raise1ex\hbox\bgroup\bgroup\Cfont}%
   \def\Rcitemark{\egroup\egroup}
   \def\LAcitemark{\def\Cfont{\ACitefont}\bgroup\ACitefont}%
   \def\RAcitemark{\egroup}
   \def\LIcitemark{\egroup\egroup}
   \def\RIcitemark{\raise1ex\hbox\bgroup\bgroup\Cfont}%
   \def\Citehyphen{\egroup--\bgroup\Cfont}
   \def\Citecomma{\egroup,\hskip0pt\bgroup%
      \Cfont}
   \def\Citebreak{}
   }

\def\Citenamedate{\Rpunct
   \def\Lcitemark{
      \def\Citebreak{\egroup\ [\bgroup\Citefont}
      \def\Citecomma{\egroup]; 
         \bgroup\let\uchyph=1\Citefont}(\bgroup\let\uchyph=1\Citefont}%
   \def\Rcitemark{\egroup])}
   \def\LAcitemark{
      \def\Citebreak{\egroup\ [\bgroup\Citefont}\def\Citecomma{\egroup], %
         \bgroup\ACitefont }\bgroup\let\uchyph=1\ACitefont}%
   \def\RAcitemark{\egroup]}
  \def\Citehyphen{\egroup--\bgroup\Citefont}
   \def\LIcitemark{\egroup}
   \def\RIcitemark{\bgroup\Citefont}
   }

%
%
%


\def\Flagstyle#1{\Flagfont#1. }

\Refstd\Citebrackets
\def\Citefont{\bf}\def\Titlefont{\sl}\def\Volfont{\bf}\def\Tomefont{\Reffont}

\def\journalarticleformat{\Reffont\let\uchyph=1\parindent=1.25pc\def\Comma{}%

\sfcode`\.=1000\sfcode`\?=1000\sfcode`\!=1000\sfcode`\:=1000\sfcode`\;=1000\sfcode`\,=1000
                \par\vfil\penalty-200\vfilneg
      \if\Ftest\present\Flagstyle\Fstr\fi%
       \if\Atest\present\bgroup\Authfont\Astr\egroup\def\Comma{\unskip, }\fi%
        \if\Ttest\present\Comma\bgroup\Titlefont\Tstr\egroup\def\Comma{, }\fi%
         \if\etest\present\hskip.2em(\bgroup\estr\egroup)\def\Comma{\unskip,
}\fi%
          \if\Jtest\present\Comma\bgroup\Tomefont\Jstr\/\egroup\def\Comma{,
}\fi%

\if\Vtest\present\if\Jtest\present\hskip.2em\else\Comma\fi\bgroup\Volfont\Vstr\egroup\def\Comma{, }\fi%
            \if\Dtest\present\hskip.2em(\bgroup\Dstr\egroup)\def\Comma{, }\fi%
             \if\Ptest\present\bgroup, \Pstr\egroup\def\Comma{, }\fi%
              \if\ttest\present\Comma\bgroup\Titlefont\tstr\egroup\def\Comma{,
}\fi%

\if\jtest\present\Comma\bgroup\Tomefont\jstr\/\egroup\def\Comma{, }\fi%

\if\vtest,\present\if\jtest\present\hskip.2em\else\Comma\fi\bgroup\Volfont\vstr\egroup\def\Comma{, }\fi%
                 \if\dtest\present\hskip.2em(\bgroup\dstr\egroup)\def\Comma{,
}\fi%
                  \if\ptest\present\bgroup, \pstr\egroup\def\Comma{, }\fi%
                   \if\Gtest\present{\Comma Gov't ordering no.
}\bgroup\Gstr\egroup\def\Comma{, }\fi%
                    \if\Mtest\present\Comma MR
\#\bgroup\Mstr\egroup\def\Comma{, }\fi%
                     \if\Otest\present{\Comma\bgroup\Ostr\egroup.}\else{.}\fi%
                      \vskip3ptplus1ptminus1pt}

\def\conferencereportformat{\Reffont\let\uchyph=1\parindent=1.25pc\def\Comma{}%

\sfcode`\.=1000\sfcode`\?=1000\sfcode`\!=1000\sfcode`\:=1000\sfcode`\;=1000\sfcode`\,=1000
                \par\vfil\penalty-200\vfilneg
      \if\Ftest\present\Flagstyle\Fstr\fi%
       \if\Atest\present\bgroup\Authfont\Astr\egroup\def\Comma{\unskip, }\fi%
        \if\Ttest\present\Comma\bgroup\Titlefont\Tstr\egroup\def\Comma{, }\fi%
         \if\Jtest\present\Comma\bgroup\Tomefont\Jstr\/\egroup\def\Comma{,
}\fi%
          \if\Ctest\present\Comma\bgroup\Cstr\egroup\def\Comma{, }\fi%
           \if\Dtest\present\hskip.2em(\bgroup\Dstr\egroup)\def\Comma{, }\fi%
            \if\Mtest\present\Comma MR \#\bgroup\Mstr\egroup\def\Comma{, }\fi%
             \if\Otest\present{\Comma\bgroup\Ostr\egroup.}\else{.}\fi%
              \vskip3ptplus1ptminus1pt}

\def\bookarticleformat{\Reffont\let\uchyph=1\parindent=1.25pc\def\Comma{}%

\sfcode`\.=1000\sfcode`\?=1000\sfcode`\!=1000\sfcode`\:=1000\sfcode`\;=1000\sfcode`\,=1000
                \par\vfil\penalty-200\vfilneg
      \if\Ftest\present\Flagstyle\Fstr\fi%
       \if\Atest\present\bgroup\Authfont\Astr\egroup\def\Comma{\unskip, }\fi%
        \if\Ttest\present\Comma\bgroup\Titlefont\Tstr\egroup\def\Comma{, }\fi%
         \if\etest\present\hskip.2em(\bgroup\estr\egroup)\def\Comma{\unskip,
}\fi%
          \if\Btest\present\Comma in
\bgroup\Tomefont\Bstr\/\egroup\def\Comma{\unskip, }\fi%
           \if\otest\present\ \bgroup\ostr\egroup\def\Comma{, }\fi%
            \if\Etest\present\Comma\bgroup\Estr\egroup\unskip,
\ifnum\Ecnt>1eds.\else ed.\fi\def\Comma{, }\fi%
             \if\Stest\present\Comma\bgroup\Sstr\egroup\def\Comma{, }\fi%
              \if\Vtest\present\Comma vol. \bgroup\Vstr\egroup\def\Comma{,
}\fi%
               \if\Ntest\present\Comma no. \bgroup\Nstr\egroup\def\Comma{,
}\fi%
                \if\Itest\present\Comma\bgroup\Istr\egroup\def\Comma{, }\fi%
                 \if\Ctest\present\Comma\bgroup\Cstr\egroup\def\Comma{, }\fi%
                  \if\Dtest\present\Comma\bgroup\Dstr\egroup\def\Comma{, }\fi%
                   \if\Ptest\present\Comma\Pstr\def\Comma{, }\fi%

\if\ttest\present\Comma\bgroup\Titlefont\Tstr\egroup\def\Comma{, }\fi%
                     \if\btest\present\Comma in
\bgroup\Tomefont\bstr\egroup\def\Comma{, }\fi%
                       \if\atest\present\Comma\bgroup\astr\egroup\unskip,
\if\acnt\present eds.\else ed.\fi\def\Comma{, }\fi%
                        \if\stest\present\Comma\bgroup\sstr\egroup\def\Comma{,
}\fi%
                         \if\vtest\present\Comma vol.
\bgroup\vstr\egroup\def\Comma{, }\fi%
                          \if\ntest\present\Comma no.
\bgroup\nstr\egroup\def\Comma{, }\fi%

\if\itest\present\Comma\bgroup\istr\egroup\def\Comma{, }\fi%

\if\ctest\present\Comma\bgroup\cstr\egroup\def\Comma{, }\fi%

\if\dtest\present\Comma\bgroup\dstr\egroup\def\Comma{, }\fi%
                              \if\ptest\present\Comma\pstr\def\Comma{, }\fi%
                               \if\Gtest\present{\Comma Gov't ordering no.
}\bgroup\Gstr\egroup\def\Comma{, }\fi%
                                \if\Mtest\present\Comma MR
\#\bgroup\Mstr\egroup\def\Comma{, }\fi%

\if\Otest\present{\Comma\bgroup\Ostr\egroup.}\else{.}\fi%
                                  \vskip3ptplus1ptminus1pt}

\def\bookformat{\Reffont\let\uchyph=1\parindent=1.25pc\def\Comma{}%

\sfcode`\.=1000\sfcode`\?=1000\sfcode`\!=1000\sfcode`\:=1000\sfcode`\;=1000\sfcode`\,=1000
                \par\vfil\penalty-200\vfilneg
      \if\Ftest\present\Flagstyle\Fstr\fi%
       \if\Atest\present\bgroup\Authfont\Astr\egroup\def\Comma{\unskip, }%

\else\if\Etest\present\bgroup\def\Eand{\Aand}\def\Eandd{\Aandd}\Authfont\Estr\egroup\unskip, \ifnum\Ecnt>1eds.\else ed.\fi\def\Comma{, }%

\else\if\Itest\present\bgroup\Authfont\Istr\egroup\def\Comma{, }\fi\fi\fi%

\if\Ttest\present\Comma\bgroup\Titlefont\Tstr\/\egroup\def\Comma{\unskip, }%

\else\if\Btest\present\Comma\bgroup\Titlefont\Bstr\/\egroup\def\Comma{\unskip,
}\fi\fi%
            \if\otest\present\ \bgroup\ostr\egroup\def\Comma{, }\fi%

\if\etest\present\hskip.2em(\bgroup\estr\egroup)\def\Comma{\unskip, }\fi%
              \if\Stest\present\Comma\bgroup\Sstr\egroup\def\Comma{, }\fi%
               \if\Vtest\present\Comma vol. \bgroup\Vstr\egroup\def\Comma{,
}\fi%
                \if\Ntest\present\Comma no. \bgroup\Nstr\egroup\def\Comma{,
}\fi%
                 \if\Atest\present\if\Itest\present
                         \Comma\bgroup\Istr\egroup\def\Comma{\unskip, }\fi%
                      \else\if\Etest\present\if\Itest\present
                              \Comma\bgroup\Istr\egroup\def\Comma{\unskip,
}\fi\fi\fi%
                     \if\Ctest\present\Comma\bgroup\Cstr\egroup\def\Comma{,
}\fi%
                      \if\Dtest\present\Comma\bgroup\Dstr\egroup\def\Comma{,
}\fi%

\if\ttest\present\Comma\bgroup\Titlefont\tstr\egroup\def\Comma{, }%

\else\if\btest\present\Comma\bgroup\Titlefont\bstr\egroup\def\Comma{, }\fi\fi%

\if\stest\present\Comma\bgroup\sstr\egroup\def\Comma{, }\fi%
                           \if\vtest\present\Comma vol.
\bgroup\vstr\egroup\def\Comma{, }\fi%
                            \if\ntest\present\Comma no.
\bgroup\nstr\egroup\def\Comma{, }\fi%

\if\itest\present\Comma\bgroup\istr\egroup\def\Comma{, }\fi%

\if\ctest\present\Comma\bgroup\cstr\egroup\def\Comma{, }\fi%

\if\dtest\present\Comma\bgroup\dstr\egroup\def\Comma{, }\fi%
                                \if\Gtest\present{\Comma Gov't ordering no.
}\bgroup\Gstr\egroup\def\Comma{, }\fi%
                                 \if\Mtest\present\Comma MR
\#\bgroup\Mstr\egroup\def\Comma{, }\fi%

\if\Otest\present{\Comma\bgroup\Ostr\egroup.}\else{.}\fi%
                                   \vskip3ptplus1ptminus1pt}

\def\technicalreportformat{\Reffont\let\uchyph=1\parindent=1.25pc\def\Comma{}%

\sfcode`\.=1000\sfcode`\?=1000\sfcode`\!=1000\sfcode`\:=1000\sfcode`\;=1000\sfcode`\,=1000
                \par\vfil\penalty-200\vfilneg
      \if\Ftest\present\Flagstyle\Fstr\fi%
       \if\Atest\present\bgroup\Authfont\Astr\egroup\def\Comma{\unskip, }%

\else\if\Etest\present\bgroup\def\Eand{\Aand}\def\Eandd{\Aandd}\Authfont\Estr\egroup\unskip, \ifnum\Ecnt>1eds.\else ed.\fi\def\Comma{, }%

\else\if\Itest\present\bgroup\Authfont\Istr\egroup\def\Comma{, }\fi\fi\fi%
          \if\Ttest\present\Comma\bgroup\Titlefont\Tstr\egroup\def\Comma{,
}\fi%
           \if\Atest\present\if\Itest\present
                   \Comma\bgroup\Istr\egroup\def\Comma{, }\fi%
                \else\if\Etest\present\if\Itest\present
                        \Comma\bgroup\Istr\egroup\def\Comma{, }\fi\fi\fi%
            \if\Rtest\present\Comma\bgroup\Rstr\egroup\def\Comma{, }\fi%
             \if\Ctest\present\Comma\bgroup\Cstr\egroup\def\Comma{, }\fi%
              \if\Dtest\present\Comma\bgroup\Dstr\egroup\def\Comma{, }\fi%
               \if\ttest\present\Comma\bgroup\Titlefont\tstr\egroup\def\Comma{,
}\fi%
                \if\itest\present\Comma\bgroup\istr\egroup\def\Comma{, }\fi%
                 \if\rtest\present\Comma\bgroup\rstr\egroup\def\Comma{, }\fi%
                  \if\ctest\present\Comma\bgroup\cstr\egroup\def\Comma{, }\fi%
                   \if\dtest\present\Comma\bgroup\dstr\egroup\def\Comma{, }\fi%
                    \if\Gtest\present{\Comma Gov't ordering no.
}\bgroup\Gstr\egroup\def\Comma{, }\fi%
                     \if\Mtest\present\Comma MR
\#\bgroup\Mstr\egroup\def\Comma{, }\fi%
                      \if\Otest\present{\Comma\bgroup\Ostr\egroup.}\else{.}\fi%
                       \vskip3ptplus1ptminus1pt}

\def\otherformat{\Reffont\let\uchyph=1\parindent=1.25pc\def\Comma{}%

\sfcode`\.=1000\sfcode`\?=1000\sfcode`\!=1000\sfcode`\:=1000\sfcode`\;=1000\sfcode`\,=1000
                \par\vfil\penalty-200\vfilneg
      \if\Ftest\present\Flagstyle\Fstr\fi%
       \if\Atest\present\bgroup\Authfont\Astr\egroup\def\Comma{\unskip, }%

\else\if\Etest\present\bgroup\def\Eand{\Aand}\def\Eandd{\Aandd}\Authfont\Estr\egroup\unskip, \ifnum\Ecnt>1eds.\else ed.\fi\def\Comma{, }%

\else\if\Itest\present\bgroup\Authfont\Istr\egroup\def\Comma{, }\fi\fi\fi%
          \if\Ttest\present\Comma\bgroup\Titlefont\Tstr\egroup\def\Comma{,
}\fi%
            \if\Atest\present\if\Itest\present
                    \Comma\bgroup\Istr\egroup\def\Comma{, }\fi%
                 \else\if\Etest\present\if\Itest\present
                         \Comma\bgroup\Istr\egroup\def\Comma{, }\fi\fi\fi%
                 \if\Ctest\present\Comma\bgroup\Cstr\egroup\def\Comma{, }\fi%
                  \if\Dtest\present\Comma\bgroup\Dstr\egroup\def\Comma{, }\fi%
                   \if\Gtest\present{\Comma Gov't ordering no.
}\bgroup\Gstr\egroup\def\Comma{, }\fi%
                    \if\Mtest\present\Comma MR
\#\bgroup\Mstr\egroup\def\Comma{, }\fi%
                     \if\Otest\present{\Comma\bgroup\Ostr\egroup.}\else{.}\fi%
                      \vskip3ptplus1ptminus1pt}

\message {)}\message {DOCUMENT TEXT}
\magnification=\magstep1
\NoBlackBoxes
\loadeufm
\baselineskip=18pt
\parskip=6pt
\def\bG{\bold G}
\def\bT{\bold T}
\def\bB{\bold B}
\def\bU{\bold U}
\def\bA{\bold A}
\def\bP{\bold P}
\def\bM{\bold M}
\def\bN{\bold N}
\def\bZ{\bold Z}

\def\bR{\Bbb R}
\def\BbbC{\Bbb C}
\def\BbbZ{\Bbb Z}
\def\calE{\Cal E}
\def\calX{\Cal X}
\topmatter
\title
$R$\snug-groups and elliptic representations for similitude groups
\endtitle
\author
David Goldberg\footnotemark"*"
\endauthor
\address
Department of Mathematics,
Purdue University,
West Lafayette, IN 47907
\endaddress
\curraddr
Mathematical Sciences Research Institute, Berkeley, CA 94720
\endcurraddr
\abstract
The tempered spectrum of the similitude groups of non-degenerate
symplectic, hermitian, or split orthogonal forms defined over
$p$\snug-adic groups of characteristic zero is studied.  The
components of representations induced from discrete series of proper
parabolic subgroups are classified in terms of $R$\snug-groups.
Multiplicity one is proved.  The tempered elliptic spectrum is
identified, and the relation between
elliptic characters appearing in a given induced representation
is determined. Those irreducible
tempered representations which are not elliptic and not fully induced from
elliptic tempered representations are described.
\endabstract
\leftheadtext{David Goldberg}
\endtopmatter
\document
\openup2\jot
\hoffset .3130 truein
\footnotetext""{*Partially supported by NSF Postdoctoral Fellowship DMS9206246}
\footnotetext""{Research at MSRI supported in part by NSF grant DMS9022140}
\footnotetext""{1991 AMS Subject Classification: 22E50, 22E35}
\noindent
\subhead Introduction \endsubhead
We continue our program of studying the explicit description of the tempered
spectrum of reductive groups defined over a $p$\snug-adic field of
characteristic zero.  The problem can be broken down into two
steps.  The first step is to classify the discrete series representations
of $G$ and of all its Levi subgroups.  The second step is to determine the
components of those representations which
are parabolically induced from a discrete series representation of a proper
parabolic subgroup $P$ of $G.$  It is this second step that we
shall concern ourselves with.  This problem
also has two distinct parts.  The first is to determine
the criteria for an induced representation to be reducible when $P$
is a maximal proper parabolic subgroup.  This
involves the computation of Plancherel measures.
The second part is to use knowledge of the rank one Plancherel
measures to construct the Knapp-Stein $R$\snug-group.  This gives one
a combinatorial algorithm for determining the structure of the induced
representations.  Those irreducible tempered representations
whose characters fail to vanish on the regular elliptic set are
called elliptic tempered representations, and are of particular interest.

In\Lspace \Lcitemark 6\Citecomma
7\Citecomma
10\Rcitemark \Rspace{}
we determined the possible $R$\snug-groups for $Sp_{2n},\ SO_n,\
U_n,$ and $SL_n.$  Here we built upon the work of Keys\Lspace \Lcitemark
18\Citecomma
19\Rcitemark \Rspace{}.
Arthur gives a criteria for determining the elliptic
tempered representations
in terms of $R$\snug-groups in\Lspace \Lcitemark 2\Rcitemark \Rspace{}.
In\Lspace \Lcitemark 15\Rcitemark \Rspace{} Herb described the elliptic
tempered representations
for $Sp_{2n}$ and $SO_n.$    For $SO_{2n}$ she
exhibited irreducible tempered representations which
are  neither elliptic nor irreducibly induced from
an elliptic tempered representation of a proper parabolic subgroup.
Such representations do not exist for $F=\R,$\Lspace \Lcitemark 21\Rcitemark
\Rspace{},
and this is the first known example of them in the $p$\snug-adic case.
Using some results of Herb and Arthur, we described the
elliptic representations for $SL_n$ and the quasi-split unitary groups
$U_n$\Lspace \Lcitemark 6\Citecomma
10\Rcitemark \Rspace{}.
The problem of determining the zeros of Plancherel measures has
been explored by Shahidi  in\Lspace \Lcitemark 24\Citecomma
25\Citecomma
26\Citecomma
27\Citecomma
28\Citecomma
29\Rcitemark \Rspace{}.  In\Lspace \Lcitemark 29\Rcitemark \Rspace{} Shahidi
determined the
reducibility criteria for the Siegel parabolic subgroups of $Sp_{2n}$ and
$SO_n.$  In\Lspace \Lcitemark 24\Rcitemark \Rspace{} he determined the
criteria for parabolics of $SO_{2n}$
whose Levi component is of the form $GL_k\times SO_{2m}$
where $k$ is even.  In\Lspace \Lcitemark 9\Rcitemark \Rspace{} we determined
the
criteria for the Siegel parabolics of $U_n.$  Our joint work with Shahidi
\Lcitemark 12\Rcitemark \Rspace{}
has determined this criteria for  symplectic and quasi-split
orthogonal groups in the case where $\M=GL_k\times\G(m),$
with $k$ even.
Our continuing joint work with Shahidi will
will examine other maximal parabolic subgroups of classical
groups.

Here we compute the possible
$R$\snug-groups for the similitude group of a non-degenerate
symplectic, symmetric, or quasi-split hermitian form, defined over $F.$
We show that all the $R$\snug-groups are elementary two groups (cf. Theorem
2.6).
Moreover, the  $R$\snug-groups are all contained in the subgroup of the Weyl
group consisting of sign changes.
For $GU_{2n}$ this extends the work of Keys\Lspace \Lcitemark 19\Rcitemark
\Rspace{}.
We then undertake the computation of the elliptic tempered representations.
Suppose $\M$ is a Levi subgroup of $\G.$
Then, for some collection of positive integers
$m_1,\dots,m_r,$ and some $m\geq 0,$ we have
$M=\M(F)\simeq GL_{m_1}\times\dots\times GL_{m_r}\times G(m),$
where
$G(m)=GSp_{2m},\ GU_{2m},\ GO_{2m+1},\ GU_{2m+1},$
or $GO_{2m}$ as appropriate (cf. Section 2).
If $\s$ is a discrete series representation
of $M,$ then the induced representation $i_{G,M}(\s)$
has elliptic constituents if and only if
the longest possible sign change $w_0=C_1\dots C_r$ is in $R(\s)$
(cf. Proposition 3.3).
We make this more explicit
on a case by case basis (cf. Theorems 3.4 and 3.6).
If $\G=GO_{2n+1}$ or $GU_{2n+1},$
then every irreducible tempered representation of $G=\G(F)$
is either elliptic, or is irreducibly induced from an elliptic tempered
representation of a proper Levi subgroup (cf. Theorem 3.4).  For $\G=GU_{2n},\
GSp_{2n},$
or $GO_{2n},$ this statement fails to be true.  We classify those irreducible
tempered representations which are non-elliptic, and are not irreducibly
induced from elliptic tempered representations
(cf. Theorems 3.4 and 3.6).

We also show that when an induced representation has elliptic components,
then the elliptic characters of those components satisfy a
particularly nice relation.  More specifically, suppose $\k$ is
an irreducible representation of the $R$\snug-group $R,$
and $\pi_\k$ is the irreducible subrepresentation of $i_{G,M}(\s)$
attached to $\k$\Lspace \Lcitemark 2\Citecomma
19\Rcitemark \Rspace{}.  Let $\T_\k^e$
be the restriction of the character of $\pi_\k$ to the regular
elliptic set.  Then
$\T_\k^e=\k(C_1\dots C_r)\T_1^e,$ where
$1$ is the trivial character (cf. Theorem 3.9).
This result is similar to the ones obtained by Herb in\Lspace \Lcitemark
15\Rcitemark \Rspace{},
and is in fact motivated by that work.

I would like to thank Freydoon Shahidi and Rebecca Herb for
their comments on this work.  I would also like to thank the
Mathematical Sciences Research Institute, in Berkeley,
for providing the pleasant atmosphere in which this work
was concluded.

\noindent
\subhead Section 1\ \ Preliminaries\endsubhead

Suppose $F$ is a nonarchimedean local field of characteristic zero and residual
characteristic $q.$
Let $\bG$ be a connected reductive quasi--split algebraic group defined over
$F.$
We denote by $G$ the $F$\snug-rational points of $\bG.$
We let $G^e$ be the collection of regular elliptic elements of $G,$ i.e., the
set of regular elements whose centralizers in $G$ are compact modulo the center
of $G.$

We denote by $\calE_c(G)$ the collection of equivalence classes of irreducible
admissible representations of $G.$
We make no distinction between an irreducible admissible representation $\pi$
of $G,$ and its class $[\pi]$ in  $\calE_c (G).$
If $\pi\in\calE_c (G),$ then the distribution character $\Theta_\pi$ of $\pi$
is given by a locally integrable function\Lspace \Lcitemark 13\Rcitemark
\Rspace{}.
We let $\Theta_\pi^e$ be the restriction of $\Theta_\pi$ to $G^e.$
We say that $\pi$ is elliptic if $\Theta_\pi^e\not= 0.$
We let $\calE_2(G),\ \calE_t (G),$ and $\calE_e (G)$ be the collection of
discrete series, tempered, and tempered elliptic classes in $\calE_c(G),$
respectively.
Then $\calE_2(G) \subset \calE_e (G)\subset \calE_t (G).$

We fix a maximal torus $\bT$ of $\bG,$ and let $\bT_d$ be the maximal split
subtorus of $\bT.$
Denote by $\Phi=\Phi(\bG,\bT_d)$ the set of roots of $\bT_d$ in $\bG.$
Let $\Delta$ be a choice of simple roots in $\Phi,$ and let
$\Phi^+=\Phi^+(\bG,\bT_d)$ be the positive roots with respect to this choice of
$\Delta.$
The choice of $\Delta$ also determines a Borel subgroup $\bB=\bT\bU$ of $\bG.$
If $\theta\subseteq \Delta,$ then we let $\bA_\theta$ be the split subtorus of
$\bT_d$ determined by $\theta.$
Then $\bA_\theta$ is the split component of a unique parabolic subgroup
$\bP_\theta=\bM_\theta \bN_\theta$ of $\bG$ containing $\bB.$

Suppose $\bA=\bA_\theta$ for some $\theta.$
Then $\bM=\bM_\theta=Z_\bG (\bA)$ is a Levi subgroup of $\bG$ with split
component $\bA.$
We denote by $W(\bG,\bA)$ the Weyl group $N_\bG (\bA)/\bM.$
If $w\in W(\bG,\bA),$ we let $\tilde{w}\in N_\bG (\bA)$ be a representative for
$w.$

Suppose $(\sigma,V)\in \calE_2 (M).$
Let
$$
V(\sigma)=\{ f\in C^\infty (G,V)|f(mng)=\sigma(m) \delta_P^{1/2} (m) f(g),\
\forall m\in M,\ n\in N,\ g\in G\}.
$$
Here $\d_P$ denotes the modular
function of $P.$
The unitarily induced representation, $\Ind_P^G(\sigma),$ of $G$ on $V(\sigma)$
is given by right translation.
Since the class of $\Ind_P^G(\sigma)$ depends only on $\bM$ and not on the
choice of $\bN,$ we may denote its class by $i_{G,M}(\sigma).$
We now recall the theory of $R$\snug-groups which allows us to determine the
structure of $\Ind_P^G (\sigma).$
For more details see\Lspace \Lcitemark 7\Citecomma
11\Citecomma
18\Citecomma
25\Citecomma
30\Citecomma
31\Rcitemark \Rspace{}.
For $\sigma\in\calE_2 (M),$ we let $W(\sigma)=\{w\in
W(\G,\A)|\tilde{w}\sigma=\sigma\}.$
Here $\tilde{w}\sigma(m)=\sigma(\tilde{w}^{-1} m \tilde{w}).$
Since the class of $\tilde{w}\sigma$ is independent of the choice of
$\tilde{w},$ we may denote it by $w\sigma.$
Let $\Phi(\bP,\bA)$ be the reduced roots of $\bA$ in $\bP.$
For $\beta\in \Phi(\bP,\bA),$ we let $\bA_\beta$ be the subtorus of $\bA$
defined by $\beta.$
Let $\bM_\beta=Z_\bG (\A_\beta).$
Then $\bM$ is a maximal Levi subgroup of $\bM_\beta.$
Let $\bN_\beta=\bM_\beta \cap \bN,$ and $^*\bP_\beta=\bM \bN_\beta.$
Then $^*\bP_\beta$ is a maximal parabolic subgroup of $\bM_\beta.$
Let $\mu_\beta(\sigma)$ be the Plancherel measure of $\sigma$ with respect to
$\beta$\Lspace \Lcitemark 14\Citecomma
31\Rcitemark \Rspace{}.
The value of $\mu_\beta$ is given by ratios of certain values of Langlands
$L$\snug-functions\Lspace \Lcitemark 28\Rcitemark \Rspace{}.
For our purposes, it is enough to know that $\mu_\beta(\sigma)=0$ if and only
if $W(\bM_\beta,\bA) \cap W(\sigma)\not= \{1\}$ and $i_{M_\beta,M}(\sigma)$ is
irreducible.
Let $\Delta'=\{\beta\in\Phi (\bP,\bA)|\mu_\beta(\sigma)=0\}.$
Then $\pm \Delta'$ is a subroot system of $\Phi$\Lspace \Lcitemark
5\LIcitemark{}, VI, \Sec2, Proposition 9\RIcitemark \Rcitemark \Rspace{}.
Thus, the group $W'$ generated by the reflections determined by the elements of
$\Delta'$ is a subgroup of $W(\sigma).$
Let $R=R(\sigma)=\{w\in W(\sigma)|\ w\beta > 0,\ \forall \beta \in \Delta'\}.$

\proclaim{Theorem 1.1 (Knapp--Stein, Silberger\Lspace \Lcitemark 20\Citecomma
30\Citecomma
32\Rcitemark \Rspace{})}
For any $\sigma\in \calE_2 (M)$ we have $W(\sigma)=R\ltimes W'.$
Moreover, the commuting algebra $C(\sigma)$ of $\Ind_P^G(\sigma)$ has dimension
$|R|.$\qed
\endproclaim

We call $R$ the Knapp--Stein $R$\snug-group attached to $\sigma.$
Let $w\in W(\sigma),$ and choose $T_w: V @>>> V$ which defines an isomorphism
between $\sigma$ and $\tilde{w}\sigma.$
Let $w_1,w_2\in R,$ and define $\eta(w_1,w_2)$ by
$T_{w_1w_2}=\eta(w_1,w_2)T_{w_1}T_{w_2}.$
Then $\eta \colon R\times R @>>> \C^\times$ is a 2--cocycle.

\proclaim{Theorem 1.2 (Keys\Lspace \Lcitemark 19\Rcitemark \Rspace{})}
\roster
\item"(a)" The commuting algebra $C(\sigma)$ of $\Ind_P^G(\sigma)$ is
isomorphic to $\BbbC[R]_\eta,$ the group algebra of $R$ twisted by the cocycle
$\eta.$

\item"(b)" Suppose $\eta$ is a coboundary.
Let $\hat{R}$ be the set of equivalence classes of irreducible representations
of $R.$
Then there is a natural one-to--one correspondence, $\kappa \longmapsto
\pi_\kappa,$ between $\hat{R}$ and the equivalence classes of
subrepresentations of $i_{G,M}(\sigma)$ such that $\dim \kappa=\dim
\Hom\limits_G (\pi_\kappa,i_{G,M}(\sigma)).$
In particular, if $R$ is abelian and $\eta$ is a coboundary, then
$C(\sigma)\simeq \BbbC[R],$ and $i_{G,M}(\sigma)$ decomposes as $|R|$
inequivalent irreducible subrepresentations.
\qed
\endroster
\endproclaim

Let $\frak a=\Hom (X(\bM)_F,\BbbZ)$ be the real Lie algebra of $\bA.$
(Here $X(\bM)_F$ is the collection of $F$\snug-rational characters of $\bM$.)
For $w\in W(\G,\A),$ we let $\frak a_w=\{H\in\frak a|w\cdot H=H\}.$
Let $\bZ$ be the split component of $\bG,$ and let $\frak z$ be its real Lie
algebra.
Suppose $\sigma\in \calE_2(M),$ and let $R$ be the $R$\snug-group of $\sigma.$
We let $\frak a_R=\bigcap\limits_{w\in R}\frak a_w.$
In\Lspace \Lcitemark 2\Rcitemark \Rspace{} Arthur gives an explicit criteria
determining
when $i_{G,M}(\s)$ has an elliptic subrepresentation.
We give a weak version of his results which is sufficient for our needs.

\proclaim{Theorem 1.3 (Arthur \Lspace \Lcitemark 2\Rcitemark \Rspace{})}
Suppose $R$ is abelian and $\eta$ is a coboundary.
Then the following are equivalent:

\roster
\item"(a)"$i_{G,M}(\sigma)$ has an elliptic constituent;

\item"(b)"all the constituents of $i_{G,M}(\sigma)$ are elliptic;

\item"(c)"there is a $w\in R$ with $\frak a_w=\frak z.$\qquad \qed
\endroster
\endproclaim

We are also interested in which tempered representations are not elliptic, and
whether they appear as irreducibly induced from tempered representations.
If a tempered representation is irreducibly induced, we wish to
know when the inducing representation is elliptic.
The next two results allow us to determine
these things.  We again state
a weaker version of Arthur's result
(Theorem 1.4) than the one given in\Lspace \Lcitemark 2\Rcitemark \Rspace{}.

Suppose $\M'\supset \M$ is a Levi subgroup of $G$
satisfying Arthur's
compatibility condition
with respect to $\DE'$\Lspace \Lcitemark 2\LIcitemark{}, \Sec 2\RIcitemark
\Rcitemark \Rspace{}.
Then $R'=R\cap  W(\M',\A)$ is the $R$\snug-group attached to
$i_{M',M}(\s).$
If
$\kappa'\in\Hat{R'},$ then we let $\tau_{\kappa'}$ to be the irreducible
component of $i_{M',M}(\sigma)$ corresponding to $\kappa'.$  Let
$$\hat R(\kappa')=\{\kappa\in\hat R\ |\ \kappa(w)=\kappa'(w),\forall w\in
R'\}.$$
We let $\pi_\kappa$ be the irreducible constituent of $i_{G,M}(\sigma)$
corresponding to $\kappa.$

\proclaim{Theorem 1.4 (Arthur\Lspace \Lcitemark 2\LIcitemark{},
\Sec2\RIcitemark \Rcitemark \Rspace{})}
Suppose that $R$ is abelian and $C(\s)\simeq\C[R].$ Then,
for any
$\kappa'\in\Hat{R'},$
we have
$$i_{G,M'}(\tau_{\kappa'})=\underset{\kappa\in\hat R(\kappa')}\to\bigoplus
\pi_\kappa. \qed$$
\endproclaim

\proclaim{Proposition 1.5 (Herb \Lspace \Lcitemark 15\Rcitemark \Rspace{})}
Suppose $R$ is abelian and $C(\sigma)\simeq
{\Bbb C}[R].$  Let $\pi$ be an irreducible constituent of $i_{G,M}(\sigma).$
Then  $\pi=i_{G,M'}(\tau)$ for
a proper Levi subgroup $\M'$ and some $\tau\in\Cal E_t(M'),$ if
and only if $\frak a_R\not=\z.$
Moreover,
$\M'$ and $\tau$ can be chosen so that $\tau$ is elliptic if and only if
there is a $w_0\in R$ with $\frak a_R=\frak a_{w_0}.$\qed
\endproclaim

\subhead Section 2.\ \ The Similitude Groups and their
$R$\snug-groups\endsubhead

We now describe the possible $R$\snug-groups that can arise when $\bG$ is one
of the following groups: $GSp_{2n},\ GO_n,$ or $GU_n.$
Thus, $\bG$ is the similitude group of a non--degenerate symplectic, symmetric,
or hermitian form.
In the last case we assume that the hermitian form defines the quasi--split
unitary group of rank $[{n\over 2}],$ and we let $E/F$ be the quadratic
extension of $F$ over which $\bG$ splits.
Denote by $x \mapsto \overline{x}$ the Galois automorphism of $E$ over $F.$
We let $N(x)=x\overline{x}$ be the norm map.

Let $J$ be the form with respect to which $\bG$ is defined.
Then, if $\bG=GSp_{2n}$ or $GO_n,$
$$
\bG=\{g\in GL_n|\ ^t\!g Jg=\lambda(g) J\text{ for some }\lambda(g)\in \bG_m\}.
$$
For the unitary similitude groups
$$
\bG=\{g\in \Res\nolimits_{E/F} GL_n|\ ^t\overline{g} Jg=\lambda(g) J,\text{ for
some }\lambda(g)\in \Res\nolimits_{E/F} \bG_m\}.
$$
In each case we call $\lambda$ the multiplier character of $\bG.$
If $\o$ is a character of $F^\times,$ then we also denote by $\o$
the character of $\G(F)$ given by $\o\circ\l.$
We fix the following forms for the various $\bG.$
We let
$$
w_\ell=\pmatrix &&&&&&1\\
                &&&&&\cdot\\
                &&&&\cdot\\
                &&&\cdot\\
                &&1\\
                &1\\
                1\endpmatrix \in M_\ell (F).
$$
Then we take $J=\pmatrix 0&w_n\\
-w_n&0\endpmatrix$ if $\bG=GSp_{2n}.$
For $\bG=GO_n,$ we take $J=w_n.$
For the unitary groups, we fix an element $y\in E$ with $\overline{y}=-y.$
Then we let
$$
\align
J&=y\pmatrix 0&w_n\\ -w_n&0\endpmatrix \text{ if } \bG=GU_{2n}\text{ and}\\
J&=\pmatrix 0&0&y w_n\\
            0&1&0\\
            y(-w_n)&0&0\endpmatrix \text{ if }\bG=GU_{2n+1}.
\endalign
$$

We will often denote $GSp_{2n},\ GO_{2n},\ GO_{2n+1},\ GU_{2n},$ or $GO_{2n+1}$
by $\bG(n).$
We formally let $\bG(0)=GL_1$ if $\bG\neq GU_{2n+1}$
and  $\bG(0)=\Res\nolimits_{E/F} GL_1$ if $\bG=GU_{2n+1}.$

Let $\bT$ be the maximal torus of diagonal elements in $\bG.$
We given an explicit description of $T=\bT(F)$ in each case.
We adopt the convention of denoting a diagonal matrix by
$\diag\{(x_1,\ldots,x_k)\}.$
If $\bG=GSp_{2n}$ or $GO_{2n},$ then
$$
T=\{\diag\{x_1,\ldots,x_n,\lambda x_n^{-1},\ldots,\lambda
x_1^{-1}\}|x_i,\lambda\in F^\times\}.\tag2.1
$$
If $\bG=GO_{2n+1},$ then
$$
T=\{\diag\{x_1,\ldots,x_n,\lambda,\lambda^2 x_n^{-1},\ldots,\lambda^2
x_1^{-1}\}|x_i,\lambda\in F^\times\}\tag2.2.
$$
If $\bG=GU_{2n},$ then
$$
T=\{\diag(x_1,\ldots,x_n,\lambda \overline{x}_n^{-1}, \ldots, \lambda
\overline{x}_1^{-1})|\lambda\in F^\times, x_i\in E^\times\},\tag2.3
$$
While if $\bG=GU_{2n+1},$ then
$$
T=\{\diag\{x_1,\ldots,x_n,\lambda,N(\lambda)\overline{x}_n^{-1},\ldots,N(\lambda)\overline{x}_1^{-1}\}|x_i,\lambda \in E^\times\}.\tag2.4
$$
In each case we let $\bT_d$ be the maximal $F$ split subtorus of $\bT.$
For $\bG=GSp_{2n}$ or $GO_n,$ we have $\bT=\bT_d.$
For $\bG=GU_{2n},$ $T_d$ is given by (2.1), while $T_d$ for $GU_{2n+1}$ is
given by (2.2).
We sometimes denote an element of $T$ by $t(x_1,\ldots,x_n,\lambda).$

Note that in each case $W(\bG,\bT_d)$ is given by $S_n \ltimes \BbbZ^n_2.$
We use standard cycle notation for the elements of $S_n.$
Thus,
$$
(ij)\colon t(x_1,\ldots,x_n,\lambda)\longmapsto
t(x_1,\ldots,x_j,\ldots,x_i,\ldots,x_n,\lambda).
$$
We let $c_i$ be the $i$\snug-th sign change,
$$
c_i\colon t(x_1,\ldots,x_n,\lambda)\longmapsto t(x_1,\ldots,x_{i-1},\lambda
x_i^{-1},\ldots,x_n,\lambda).
$$
We call any product of the $c_i$'s a sign change.
We remark that for $\bG=SO_{2n},$ only those sign changes which are products of
an even number of $c_i$'s are in $W(\bG,\bT_d).$
However, if
$\bG=GO_{2n},$ then the sign changes $c_i$ can be represented in $\bG.$
Finally we remark that $W(\bG,\bT_d)=\langle (ij)\rangle \ltimes \langle
c_i\rangle.$

Let $\Delta$ be the choice of simple roots of $\bT_d$ in $\bG$ for which the
associated Borel subgroup consists of upper triangle matrices.
Suppose $\theta \subset \Delta.$
Let $\bP=\bP_\theta$ be the associated parabolic subgroup, and write
$\bP=\bM\bN$ for the Levi decomposition of $\bP.$

For some positive integers $m_1,\ldots,m_r,$ and some $m\geq 0$ we have,
$$
M=\bM(F) \simeq GL_{m_1}\times \ldots \times GL_{m_r} \times G(m),\tag2.5
$$
where
$$
GL_{m_i}=\cases GL_{m_i}(F),&\text{if $\bG\not= GU_n$};\\
GL_{m_i}(E),&\text{if $\bG=GU_n$}.\endcases
$$
In particular, we may assume that $M$ consists of block diagonal matrices.
Let $\varepsilon\colon GL_{m_i} @>>> GL_{m_i}$ be given by
$\e(g)=\ ^t\overline{g}^{-1}$ (where for $GL_{m_i}(F),$ the
Galois automorphism is of course trivial), and we let
$\varepsilon'(g)=w_{m_i}\varepsilon(g)w_{m_i}^{-1}.$
Then we may assume that
$$
M=\left\{\diag\{g_1,\ldots,g_r,g,\lambda(g)\e'(g_r),\ldots,\lambda(g)
\varepsilon'(g_1)\}\Big|\ g_i\in GL_{m_i},g\in G(m)\right\}.
$$
We may denote an element of $M$ by $(g_1,\ldots,g_r,g).$
Note that
$$
A=\bA(F)=\{\diag\{x_1 I_{m_1},\ldots,x_r I_{m_r}, I_k, \overline{x}_r^{-1}
I_{m_r} \ldots \overline{x}_1^{-1} I_{m_1}\}\},
$$
where $k$ is chosen appropriately to be $2m$ or $2m+1.$
For $1\leq i\leq r,$ we let $b_i=\sum\limits_{j=1}^i m_j,$
and set $b_0=0.$
For $1\leq i\leq j\leq r-1,$ we let $\alpha_{ij}=e_{b_i}-e_{b_j+1},$ and
$\beta_{ij}=e_{b_i}+e_{b_j+1}.$
For $1\leq i\leq r,$ we set
$$
\gamma_i=\cases 2e_{b_i},&\text{if $G=GSp_{2n}$ or $GU_{2n}$;}\\
e_{b_i},&\text{if $G=GO_{2n+1}$ or $GU_{2n+1}$;}\\
e_{b_i-1}+e_{b_i},&\text{if $G=GO_{2n}$}.\endcases
$$
Then $\Phi(\bP,\bA)=\{\alpha_{ij},\beta_{ij}\}_{1\leq i\leq j\leq r-1}\cup
\{\gamma_i\}_{1\leq i\leq r}.$

Note that the Weyl group $W(\bG,\bA)$ is a subgroup of $S_r \ltimes \BbbZ_2^r.$
Namely, if $m_i=m_j,$ then the permutation:
$$
w_{ij}=\prod_{k=1}^{m_i} (b_{i-1}+k\ b_{j-1}+k)
$$
is in $W(\bG,\bA),$ and the map $w_{ij} \mapsto (ij)$ gives an isomorphism of
this permutation group with a subgroup of $S_r.$
For each $i,$ we have the sign change
$$
C_i=\prod_{j=1}^{m_i} c_{b_{i-1}+j}
$$
is in $W(\bG,\bA).$
We call $C_i$ a block sign change.
The $C_i$'s generate the subgroup $\BbbZ_2^r$ in the semidirect product.
Note that, for $(g_1,\ldots,g_r,g)\in M,$ we have
$$
\align
w_{ij} (g_1,\ldots g_r,g)&=(g_1,\ldots g_j,\ldots g_i \ldots g_r,g),\text{
and}\\
C_i(g_1 \ldots g_r,g)&=(g_1,\ldots\lambda(g)\varepsilon(g_i),\ldots,g_r,g).
\endalign
$$
Suppose $\sigma\in\calE_2 (M).$
Then $\sigma\simeq \sigma_1\otimes \sigma_2\otimes \ldots\otimes \sigma_r
\otimes\rho,$ for some $\sigma_i\in \calE_2 (GL_{m_i})$ and $\rho\in\calE_2
(G(m)).$
Let $\omega_i$ be the central character of $\sigma_i.$
Note that if $w_{ij}\in W(\bG,\bA),$ then
$$
w_{ij}\sigma\simeq \sigma_1\otimes\ldots \otimes \sigma_j \otimes \ldots
\otimes \sigma_i \otimes \ldots \otimes \sigma_r \otimes \rho.
$$
If $1\leq i\leq r,$ then
$$
C_i\sigma \simeq \sigma_1 \otimes \ldots \otimes \sigma_i^\varepsilon \otimes
\ldots \otimes \sigma_r \otimes (\rho\otimes \omega_i).
$$
Here $\s_i^\e(g_i)=\s_i(\e(g_i)).$
Thus, $W(\sigma)\not= \{1\}$ if and only if at least one of the following
holds:
$$
\align
\sigma_i &\simeq \sigma_j \text{ for some } i\not= j,\tag2.6\\
\sigma_i &\simeq \sigma_j^\varepsilon\text{ and }\rho\otimes\o_i\o_j\simeq\rho,
\text{ or },\tag2.7
\endalign
$$
for some subset $\{i_j\}\subseteq \{1,\ldots,r\},$
$$
\sigma_{i_j}\simeq \sigma_{i_j}^\varepsilon \text{ for all } j\text{ and }
\rho\otimes\left(\prod_j \omega_{i_j}\right)\simeq\rho.\tag2.8
$$
Note that if (2.6) holds then $w_{ij}\in W(\sigma).$
If (2.7) holds then $w_{ij} C_i C_j \in W(\sigma).$
Finally, (2.8) implies $\prod\limits_j C_{i_j} \in W(\sigma).$
Note that in (2.7), since $\sigma_i \simeq \sigma_j^\varepsilon,$ we have
$$\omega_i=\cases \omega_j^{-1}&\text{if }\G\neq GU_n,\\ \bar\o_j^{-1}&\text{if
}\G=GU_n.\endcases$$
Since
$$
\lambda(G(m))=\cases F^\times,&\text{if }\bG=GSp_{2n},\ GO_{2n},\text{ or
}GU_{2n};\\
F^{x^2},&\text{if }\bG=GO_{2n+1};\\
N_{E/F}E^\times,&\text{if }\bG=GU_{2n+1},\endcases
$$
(see\Lspace \Lcitemark 16\Citecomma
22\Rcitemark \Rspace{}),
we see that $\o_i\o_j\circ\l=1$ if $\s_i\simeq\s_j^\e.$
Thus, (2.7) simply becomes $\s_i\simeq\s_j^{\e}.$

We begin the computation of the possible $R$\snug-groups with the following
standard result.
Recall that $\Delta'=\{\beta\in\Phi(\bP,\bA)|\mu_\beta(\sigma)=0\}.$

\proclaim{Lemma 2.1} The reduced root $\alpha_{ij}\in \Delta'$ if and only if
$\sigma_i \simeq \sigma_{j+1}.$
Similarly, $\beta_{ij}\in \Delta'$ if and only if $\sigma_i \simeq
\sigma_{j+1}^\varepsilon.$
\endproclaim

\demo{Proof} Note that
$$
M_{\alpha_{ij}} \simeq \left( \prod_{k\not= i,j+1} GL_{m_k}\right)\times
GL_{m_i+m_{j+1}}\times G(m) \simeq M_{\beta_{ij}}.
$$
Thus, by the results of Ol'{\v{s}}anski{\v{i}\Lspace \Lcitemark 23\Rcitemark
\Rspace{},  Bernstein--Zelevinsky\Lspace \Lcitemark 4\Citecomma
33\Rcitemark \Rspace{}, and Jacquet\Lspace \Lcitemark 17\Rcitemark \Rspace{},
$i_{M_{\alpha_{ij}},M}(\sigma)$ and $i_{M_{\beta_{ij}},M}(\sigma)$
are both irreducible for all $\sigma\in \calE_2 (M).$
Thus, $\mu_{\alpha_{ij}}(\sigma)=0$ if and only if there is some $w\not= 1$ in
$W(\M_{\alpha_{ij}},\A)$ with $w\sigma\simeq \sigma.$
If $m_i \not= m_{j+1},$ then $W(\M_{\alpha_{ij}},\A)=\{1\}.$
However, if $m_i=m_{j+1},$ then $W(\M_{\alpha_{ij}},\A)=\{w_{ij},1\},$ and so
$\mu_{\alpha_{ij}}(\sigma)=0$ if and only if $\sigma_i\simeq \sigma_{j+1}.$
Since $\beta_{ij}=C_{j+1}\alpha_{ij},$ we see that $\mu_{\beta_{ij}}(\sigma)=0$
if and only if $\mu_{\alpha_{ij}} (C_{ij} \sigma)=0,$ which is equivalent to
$\sigma_i \simeq \sigma_{j+1}^\varepsilon.$\qed
\enddemo

We now show that the $R$\snug-group of any $\sigma\in \calE_2 (M)$ is an
elementary 2--group.
The proof of the following lemma is based on a technique of Keys\Lspace
\Lcitemark 18\Rcitemark \Rspace{}.

\proclaim{Lemma 2.2} Suppose $w=sc\in R,$ with $s\in S_r$ and $c\in \BbbZ_2^r.$
Then $s=1.$
\endproclaim

\demo{Proof} By conjugating by a sign change, we may assume that $c$ changes
the sign of at most one $e_{b_i}$ in each orbit of $s.$
Suppose $s$ has a non--trivial cycle, which we assume is $(1\ldots j+1).$
If $c$ changes no signs among $\{e_{b_1},\ldots,e_{b_{j+1}}\},$ then we see
that $w\sigma\simeq \sigma$ implies $\sigma_1\simeq \sigma_2 \simeq \ldots
\simeq \sigma_{j+1}.$
By Lemma 2.1, we see that $\alpha_{1j}\in \Delta',$ while $w\alpha_{1j} < 0.$
This contradicts the assumption that $w\in R.$
Suppose that $c$ changes the sign of $e_{b_{j+1}}.$
Then $w\sigma\simeq \sigma$ implies $\sigma_1 \simeq \ldots \simeq
\sigma_{j+1}\simeq \sigma_1^\varepsilon.$
So by Lemma 2.1, $\beta_{1j}\in \Delta'.$
But $w\beta_{1j}=s \alpha_{1j} < 0,$ and again we have a contradiction.
Thus, $s=1.$\qed
\enddemo

Suppose that, for some subset $B\subset \{1,2,\ldots,r\},$ we have
$C_B=\prod\limits_{j\in B} C_j \in R.$
Set $\o_B=\prod\limits_{j\in B}\o_j.$
Then, for each $j\in B,\ \sigma_j \simeq \sigma_j^\varepsilon,$ and
$\rho\otimes\o_B\simeq\rho.$
Let
$$
N_\varepsilon F^\times=\cases N_{E/F} E^\times,&\text{if $\bG=GU_n$};\\
F^{x^2},&\text{if $\bG=GSp_{2n}$ or $GO_n$}.\endcases
$$
If $\sigma_i\simeq \sigma_i^\varepsilon$ then
$\omega_i|_{N_\varepsilon F^\times}=1.$
If $\bG=GSp_{2n},\ GO_{2n},$ or $GU_{2n},$ then $\lambda(G(m))=F^\times,$ while
if $\bG=GO_{2n+1},$ or $GU_{2n+1},$ then $\lambda(G(m))=N_\varepsilon
F^\times.$
Thus, for these last two groups, the condition
$\rho\otimes\o_B\simeq\rho$ is trivial.
Therefore, for $\bG=GO_{2n+1}$ or $GU_{2n+1},$ we see that $C_B\in R$
only if $\sigma_i\simeq \sigma_i^\e$
for each $i\in B.$

\definition{Definition 2.3} Let $\sigma\in \calE_2 (GL_k)$ and $\rho\in \calE_2
(G(m)).$
We say that condition $\calX_{m,k,G}(\sigma\otimes \rho)$ holds if
$i_{G(m+k),GL_k \times G(m)} (\sigma \otimes \rho)$ is reducible.
Note that a necessary condition for $\calX_{m,k,G}(\sigma\otimes\rho)$ to hold
is that $\rho\otimes\omega_{\sigma}\simeq\rho,$
where $\o_\s$ is the central character of $\s.$
\enddefinition

\definition{Definition 2.4} If $w\in W(\bG,\bA),$ then we let $R(w)=\{\alpha\in
\Phi (P,A)|w\alpha < 0\}.$
Note that $w\in R$ if and only if $w\in W(\sigma)$ and $R(w)\cap
\Delta'=\emptyset.$
\enddefinition

\proclaim{Lemma 2.5} Suppose $\bG=GO_{2n+1}$ or $GU_{2n+1}.$
Suppose $\bP=\bM \bN$ is a parabolic subgroup of $\bG$ with $M=\bM(F)\simeq
GL_{m_1} \times \ldots \times GL_{m_r} \times G(m).$
Let $\sigma\in \calE_2 (M),$ with $\sigma\simeq \sigma_1 \otimes \ldots \otimes
\sigma_r \otimes \rho.$
If $c\in R,$ and $c=C_B=\prod\limits_{j\in B} C_j,$ then $C_j \in R$ for each
$j\in B.$
\endproclaim

\demo{Proof} Note that, for each $j\in B,$ we have $R(C_j)\subset R(c).$
Thus, if $C_j \beta < 0$ then $C_B \beta < 0,$ and hence $\beta\notin \Delta'.$
Now it is enough to note that, since $c\sigma\simeq \sigma,$ we have $\sigma_j
\simeq \sigma_j^\varepsilon$ for each $j\in B,$ and thus $C_j \in W(\sigma).$
Therefore, $C_j\in R.$\qed
\enddemo

In order to give an explicit description of the $R$\snug-groups,
we need to define some terms.
Let $\rho\in\Cal E_2(G(m)),$ and set
$$X(\rho)=\{\chi\in\left(F^\times/N_\e
F^\times\right)^\wedge|\rho\otimes\chi\simeq\rho\}.$$
Let
$$
J_1(\sigma)=\{i|\calX_{m,m_i,G} (\sigma_i\otimes \rho) \text{ holds}\}.
$$ For each $i\in J_1(\s),$ we have
$\o_i\in X(\rho).$
For each $\chi\in (F^\times/N_\varepsilon F^\times)^\wedge\setminus X(\rho),$
let
$$J_\chi(\sigma)=\{i|\sigma_i\simeq \sigma_i^\varepsilon\text{ and }
\omega_i|_{F^\times}=\chi\}.$$
For $\chi=1$ or $\chi\not\in X(\rho),$ we set
$$I_\chi (\sigma)=\{i\in J_\chi (\sigma)|\sigma_i \not\simeq \sigma_j, \forall
i > j\},$$
and let $d_\chi=d_\chi (\sigma)=|I_\chi (\sigma)|$ be the number of
inequivalent $\sigma_i$ for $i\in J_\chi (\sigma).$

For a nonempty subset $S \subset (F^\times/N_\varepsilon F^\times)^\wedge,$ we
say that $S$ is {\bf minimally $\rho$\snug-trivial} if
$\prod\limits_{\chi \in S}\chi\in X(\rho),$ and $\prod\limits_{\chi\in S'} \chi
\not\in X(\rho)$ for any $\emptyset \subsetneq S' \subsetneq S.$
We let
$$
\align
\Lambda (\sigma)&=\\
&\{ S\subset (F^\times/N_\varepsilon F^\times)^\wedge\setminus
X(\rho)\big|S\text{ is minimally }\rho\snug-\text{trivial and }
d_\chi(\sigma)\geq 1, \forall \chi\in S\}.
\endalign
$$

Let $\{S_1,\ldots,S_k\}\subseteq \Lambda (\sigma).$
We let
$$X(\{S_1,\ldots,S_k\})=\{\chi\ |\ \chi\in S_j\text{ for an odd number of }
j\}.$$
We say that $\Lambda (\sigma)'$ is a {\bf basis} for $\Lambda(\sigma)$ if, for
every $S\in \Lambda(\sigma),$
there is some $\{S_1,\dots,S_k\}\subset \Lambda(\s)'$
with $S=X(\{S_1,\ldots,S_k\}),$ and $\Lambda(\sigma)'$ is minimal with respect
to this property.

\proclaim{Theorem 2.6} Let $\bM$ be a Levi subgroup of $\bG,$ with
$$
M\simeq GL_{m_1} \times \ldots \times GL_{m_r} \times G(m).
$$
Suppose $\sigma\in \calE_2 (M),$ with $\sigma\simeq \sigma_1 \otimes \ldots
\otimes \sigma_r \otimes \rho.$
Let $R=R(\sigma).$
\roster
\item"(a)"If $\bG=GO_{2n+1}$ or $GU_{2n+1},$ then $R\simeq \BbbZ_2^d,$ with
$d=d_1,$ i.e., $d$ is the number of inequivalent $\sigma_i$ such that
$\calX_{m,m_i,G}(\sigma \otimes \rho)$ holds.

\item"(b)"If $\bG=GSp_{2n},\ GO_{2n},$ or $GU_{2n},$ then $R\simeq \BbbZ_2^d$
with
$$
d=d_1+\sum_{\chi\not\in X(\rho)\atop d_\chi \geq 1}
(d_\chi-1)+|\Lambda(\sigma)'|,
$$
for any basis $\Lambda(\sigma)'$ of $\Lambda(\sigma).$
\endroster
\endproclaim

\demo{Proof} (a) By Lemma 2.5, is enough to show that $C_i \in R$ if and only
if $i\in I_1 (\sigma).$
Suppose $C_i\in W(\sigma).$
Note that $R(C_i)=\{\alpha_{ij},\beta_{ij}\}_{i\leq j\leq r-1}\cup
\{\gamma_i\}.$
Note that if $\sigma_{j+1}\simeq \sigma_i,$ for some $j\geq i,$ then
$\alpha_{ij}$ and $\beta_{ij}\in \Delta',$ and thus $C_i \notin R.$
Direct computation shows that
$M_{\gamma_i}\simeq (\prod\limits_{k\not= j} GL_{m_k})\times G(m+m_i).$
Since $C_i\sigma\simeq \sigma,$
we see that $\gamma_i\in \Delta'$ if and only if
$\calX_{m,m_i,G}(\sigma_i\otimes \rho)$ does not hold.
Thus, $C_i\in R$ if and only if $\calX_{m,m_i,G}(\sigma_i\otimes \rho)$ holds
and $\sigma_i \not\simeq \sigma_j$ for all $j > i,$ i.e., if and only if $i\in
I_1 (\sigma).$

(b) Suppose $B\subseteq \{1,2,\ldots,r\}$ with $C_B=\prod\limits_{j\in B}C_j
\in R.$
By the argument of part (a), we see that, for each $j\in B,\ \sigma_k\not\simeq
\sigma_j$ for all $k > j.$
Note again that $M_{\gamma_j}\simeq (\prod\limits_{k\not= j} GL_{m_k})\times
G(m+m_j).$
Therefore,
$\gamma_j\in\Delta'$ if and only if $C_j \sigma \simeq \sigma$ and
$\calX_{m,m_j,G}(\sigma\otimes\rho)$ fails.
Since $C_j\s\simeq\s$ if and only if $\s_j\simeq\s_j^\e,$
and $\omega_j\in X(\rho),$
we see that those $j\in B$ for which $\omega_j\in X(\rho)$ are all elements of
$I_1(\sigma).$
Let
$$R_1=R_1(\s)=\langle C_j\ \big|\ j\in I_1 (\sigma)\rangle.$$
Then we have just shown that $R_1(\sigma)\subseteq R(\sigma).$
If $j\in B$ and $\omega_j\not\in X(\rho),$ then $C_j \not\in W(\sigma),$ and
thus $\gamma_j \not\in \Delta'.$
Let $\chi\in (F^\times/N_\e F^\times)^\wedge\backslash X(\rho),$ and define
$$R_\chi(\sigma)=\langle C_i C_j|i,j\in I_\chi (\sigma)\rangle.$$
If $C_i C_j \in R_\chi (\sigma),$ then $\omega_i \omega_j|_{F^\times}=1,$ so
$C_i C_j \in W(\sigma).$
Note that $R(C_iC_j)=R(C_i)\cup R(C_j).$
Moreover, by the definition of $I_\chi(\sigma),$
and Lemma 2.1, we see that, for all
$k\geq i$ and $\ell\geq j,$
the reduced roots $\alpha_{ik},\ \alpha_{j\ell},\ \beta_{ik},$
and $\beta_{j\ell}$ can not be in $\Delta'.$
Since $\gamma_i,\ \gamma_j \not\in \Delta',$ we see that $C_i C_j \in R,$ and
thus $R_\chi(\sigma)\subseteq R.$

Since each $R_\chi(\s)\subset R(\s),$
we can multiply $C_B$ by an element $C'$ of
$$R_1\times\prod\limits_{\chi\not\in X(\rho)} R_\chi (\sigma)$$
so that $C_{B'}=C' C_B$ has the property that
$B'\cap I_1(\s)=\emptyset,$ and $\o_i\neq\o_j$
for all $i\neq j$ in $B'.$
Since $\o_B\in X(\rho),$
we see that
$\left\{\omega_j|_{F^\times}\ \big|j\in B'\right\}$ can be partitioned into
minimally $\rho$\snug-trivial subsets.
For each $\chi\in (N_\e F^\times\bs F^\times)^\wedge\setminus X(\rho),$
with $d_\chi>0,$
we fix an $i_\chi\in I_\chi (\sigma).$
Let $S\in \Lambda(\sigma),$ and define $C_S=\prod\limits_{\chi\in S}
C_{i_\chi}.$
Since $\prod\limits_{\chi\in S}\chi\in X(\rho),$ we see that $C_S\in
W(\sigma),$ and by the definition of $I_\chi(\sigma),$ we see that $C_S\in R.$
Suppose $C'=C_D\in R$ has the property that $\{\omega_j|_{F^\times}|j\in D\}=S$
and $|D|=|S|.$
Then
$$
C_D=C_S\cdot \prod_{j\in D} C_j\ C_{i_{\omega_j}},$$
and each $C_j C_{i_{\omega_j}}\in R_{\omega_j} (\sigma).$
Thus, all the factors in the product are in $R.$
Fix a basis $\L(\s)'$ for $\L(\s).$
Suppose $S\in \Lambda(\sigma),$ but $S\not\in \Lambda(\sigma)'.$
Then $S=X(\{S_1,\ldots,S_k\})$ for some $\{S_1,\ldots,S_k\}\subseteq
\Lambda(\sigma)'.$
By the definition of $C_S,$ and of $X(\{S_1,\ldots,S_k\}),$ we have
$C_S=C_{S_1} C_{S_2} \ldots C_{S_k}.$
Thus, for some $\{S_1,\ldots,S_\ell\}\subseteq \Lambda(\sigma)',$ and some $C''
\in \prod\limits_{\chi\not\in X(\rho)} R_\chi (\sigma),$ we have
$C_{B'}=C_{S_1}\ldots C_{S_\ell} \cdot C''.$
Therefore, if we define
$$R'(\sigma)=\langle C_S|S\in \Lambda(\sigma)'\}\rangle,$$
then we have seen that
$$
R=\left(\prod_{\chi} R_\chi(\sigma)\right)\cdot R'(\sigma).
$$
Moreover, by considering central characters, this is a direct product.
Therefore,
$$R\simeq\Z_2^{d_1}\times\left(\prod_{\chi\not\in
X(\rho)}\Z_2^{d_\chi-1}\right)\times\Z_2^{|\L(\s)'|},$$
as claimed.\qed
\enddemo

\remark{Remark} If $G=GU_{2n},$ then $N_\varepsilon F^\times=N_{E/F} E^\times$
is of index two in $F^\times,$ and so the local class field theory character
$\omega_{E/F}$ is the only non-trivial element of $(F^\times/N_\varepsilon
F^\times)^\wedge.$
Therefore,  the group $R'$ is trivial,
and $R\simeq \Z_2^d,$ with $d=d_1,$ or $d_1+d_{\omega_{E/F}}-1$ depending on
whether or not $d_{\omega_{E/F}}=0$ or not.
We note the similarity between the
$R$\snug-group structure for $GU_{2n}$
and that for $SO_{2n}.$
We will see that the theory of elliptic representations
for $GU_{2n}$ also
parallels that that of $SO_{2n}.$
\endremark

\subhead Section 3\ \ Elliptic Representations\endsubhead

We now describe the tempered elliptic representations of $G.$
We first show that the 2--cocycle $\eta$ of $R$ described in Section 1 is
trivial.
We use an idea from\Lspace \Lcitemark 26\Rcitemark \Rspace{}.

\proclaim{Lemma 3.1} Let $\bG=GSp_{2n},\ GO_n,$ or $GU_n.$
Suppose $\bM$ is a Levi subgroup of $G$ and $\sigma\in \calE_2 (M).$
Denote by $R$ the $R$\snug-group associated to $i_{G,M}(\sigma).$
Then the commuting algebra $C(\sigma)$ is isomorphic to $\C[R].$
\endproclaim

\demo{Proof} Suppose that $\G\neq GO_{2n}.$
Let $\G_1(n)=\{g\in \G(n)|\l(g)=1\}^\circ.$
That is, $\G_1(n)$ is the symplectic, special orthogonal, or quasi-split
unitary group associated to the form $J.$  Let $\M_1$
be the Levi component of $\G_1(n)$ given by $\M_1=\M\cap \G_1(n).$
If
$$\M=GL_{m_1}\times\dots\times GL_{m_r}\times \G(m),$$
then
$$\M_1=GL_{m_1}\times\dots\times GL_{m_r}\times \G_1(m),$$
Let $\PPP_1=\M_1\N.$
Suppose that $\tau$ is an irreducible subrepresentation of $\s|_{M_1}.$
Then  $\tau\simeq\s_1\otimes\dots\otimes\s_r\otimes\rho_1,$ for some
irreducible component $\rho_1$ of $\rho|_{{G_1(m)}}.$
Let $A(\nu,\tau,w)$ and $A(\nu,\s,w)$ be the standard
intertwining operators attached to either $\s$ and $\tau,$ respectively.
\Lcitemark 19\Citecomma
25\Citecomma
28\Rcitemark \Rspace{}.
Since $N\subset G_1(n),$
$(A(\nu,\s,w)f)|_{G_1}=A(\nu,\tau,w)f|_{G_1}$
for any $f\in \Ind_{P}^{G}(\tau).$
Thus, for each reduced root $\b,$ we have
$\mu_\b(\s)=\mu_\b(\tau).$
We let $\Cal A'(\tau,w)$ and $\Cal A'(\s,w)$
be the normalized self intertwining operators
(see\Lspace \Lcitemark 10\Citecomma
15\Rcitemark \Rspace{}).
Then these operators satisfy
$$\Cal A'(\s,w_1w_2)=\eta(w_1,w_2)\Cal A'(\s,w_1)\Cal A'(\s,w_2)$$
for each $w_1,w_2\in R(\s),$
and
$$\Cal A'(\tau,w_1w_2)=\eta(w_1,w_2)\Cal A'(\tau,w_1)\Cal A'(\tau,w_2)$$
for $w_1,w_2\in R(\tau).$

Note that each element of $R$ can be represented by an element
in $G_1.$ Suppose that $C_B\in R.$  Then for each $i\in B,$
$\s_i\simeq\s_i^\e.$  Thus, by the results of\Lspace \Lcitemark 7\Rcitemark
\Rspace{}
and\Lspace \Lcitemark 6\Rcitemark \Rspace{}, we have $C_B\in W(\tau).$  If
$C_B\b<0,$
then $\mu_\b(\s)\neq0,$ and thus $\mu_\b(\tau)\neq 0.$
Consequently, $C_B\in R(\tau),$ i.e., $R(\s)\subset R(\tau).$
Now let $w_1,w_2\in R(\s),$ and $f\in \Ind_{P}^{G}(\s),$
and $f_1=f|_{G_1}.$
Then
$$\eta(w_1,w_2)\Cal A'(\s,w_1)\Cal A'(\s,w_2)f|_{G_1}=\Cal
A'(w_1w_2,\s)f|{G_1}.$$
Thus,
$$\eta(w_1,w_2)\Cal A'(w_1,\tau)\Cal A'(w_2,\tau)f_1=\Cal A'(w_1w_2,\tau)f_1.$$
On the other hand, the results of\Lspace \Lcitemark 6\Citecomma
15\Rcitemark \Rspace{}
show that
$$\Cal A'(w_1,\tau)\Cal A'(w_2,\tau)=\Cal A'(w_1w_2,\tau),$$
and thus
$\eta(w_1,w_2)=1.$

We are left to prove the lemma for $\G=GO_{2n}.$  The problem with the above
argument is that there may be $R$\snug-group elements
which are not elements of the Weyl group attached to the
corresponding parabolic subgroup of $SO_{2n}.$  We can remedy this by
looking at restriction to $M_1$ as above, and inducing to $O_{2n}.$
Let $\G_1=SO_{2n},$ and take $\M_1=\M\cap \G_1,$
and $\PPP_1=\PPP\cap \G_1$ as above.
Let
$\G_2=O_{2n}.$
If $m>0,$ then for each $i$ with
$m_i$ odd, we choose the representative for $C_i$ as in\Lspace \Lcitemark
7\Rcitemark \Rspace{}.
Then $C_i$ also represents an element of $W(\G_1,\A_1),$
where this is the obvious Weyl group.  Suppose $\s\in\Cal E_2(M),$
and $\tau$ is an irreducible subrepresentation of
$\s|_{M_1}.$
In\Lspace \Lcitemark 8\Rcitemark \Rspace{} we examine the decomposition
of the representation
$$\Ind_{P_1}^{G_2}(\tau)=\Ind_{G_1}^{G_2}\left(\Ind_{P_1}^{G_1}\left(\tau\right)\right).$$

In\Lspace \Lcitemark 1\Rcitemark \Rspace{}, Arthur extends the definition of
the
unnormalized intertwining operators $A(\nu,\tau,w)$ to the
case of disconnected groups whose component group is cyclic.
Consider these operators for $\G_2.$
Since $\N\subset \G_2,$ this integral operator is again the same
one which gives the  unnormalized operators $A(\nu,\s,w)$
for $GO_{2n}.$  In the general case, Arthur did not prove the
analytic continuation of these operators.   However,
since $A(\nu,\s,w)$ can be analytically continued,
and $A(\nu,\s,w)f|_{G_2}=A(\nu,\tau,w)(f|_{G_2}),$
we see that indeed the operators for $G_2$ can be analytically continued.

Let $C_B\in R.$  Then $C_B$ can be represented in $G_2,$
and is a representative for an element of $W(\G_2,\A_1).$
Moreover, from the results of\Lspace \Lcitemark 7\Citecomma
8\Rcitemark \Rspace{},
$C_B\tau\simeq\tau.$
We claim that $C_B\in R_{G_2}(\tau),$
where this last object is the $R$\snug-group given in\Lspace \Lcitemark
8\Rcitemark \Rspace{}
which determines the structure of $\Ind_{P_1}^{G_2}(\tau).$
It is enough to show that if $i\in B,$ then $\mu_{\ga_i}(\tau)\neq 0.$
If $C_i$ represents an element of $W(\G_1,\A_1),$
then the argument used for $GO_{2n+1},\ GSp_{2n},$ and $GU_n$
above shows that $\mu_{\ga_i}(\tau)=\mu_{\ga_i}(\s).$
Thus, we are reduced to the case where $m_i$ is odd, and
$m=0.$  Then, since $W(\M_{\ga_i},\A_1)=\{1\},$
we have $\mu_{\ga_i}(\tau)\neq 0,$ as pointed out in\Lspace \Lcitemark
7\Rcitemark \Rspace{}.
Thus, $C_B\in R_{G_2}(\s).$

Now, the operators
$\Cal A'(w,\s)$ restrict to elements $\Cal A'(\tau,w)$ of the commuting algebra
of $\Ind_{P_1}^{G_2}(\tau).$  Thus, the cocycle $\eta$
also determines the composition of these operators.
 However, these
operators can also be defined using operators which
intertwine $\tau$ and $w\tau,$ and $\eta$
is also the cocycle that arises in this way.
The argument given in Proposition 2.3 of\Lspace \Lcitemark 15\Rcitemark
\Rspace{} shows that,
we can choose $T_w'$
for each $w$ in $R_{G_2}(\tau),$ satisfying
$T_w'\tau=w\tau T_w',$
with $T_{w_1w_2}'=T_{w_1}'T_{w_2}'.$
But this shows that
$$\Cal A'(w_1w_2,\tau)=\Cal A'(w_1,\tau)\Cal A'(w_2,\tau),$$
and therefore $\eta\equiv 1.$
This completes the lemma.
\qed
\enddemo

\proclaim{Corollary 3.2} For any $\sigma\in \calE_2 (M),$ the induced
representation $i_{G,M}(\sigma)$
decomposes as a direct sum of $|R|$ inequivalent representations.\qed
\endproclaim

\proclaim{Proposition 3.3} Let $\bG=GSp_{2n},\ GO_n,$ or $GU_n.$
Suppose $\bM$ is a Levi subgroup of $\bG,$ with $M\simeq GL_{m_1}\times \ldots
\times GL_{m_r} \times G(m).$
Let $\sigma\in \calE_2 (M)$ with $\sigma\simeq \sigma_1 \otimes \ldots \otimes
\sigma_r \otimes \rho.$
Suppose $R$ is the $R$\snug-group attached to $i_{G,M}(\sigma).$
Then the following are equivalent:
\roster
\item"(a)"$i_{G,M}(\sigma)$ has an elliptic constituent;

\item"(b)"all the constituent of $i_{G,M}(\sigma)$ are elliptic;

\item"(c)"$C_1\ldots C_r \in R.$
\endroster
\endproclaim

\demo{Proof} By Lemma 3.1, and Theorems 1.3 and 2.6,
conditions (a) and (b) are equivalent.
Furthermore, both are equivalent to $\frak a_w=\frak z=\{0\}$ for some $w\in
R.$
Note that
$$
\frak a=\{\diag\{x_1 I_{m_1},\ldots,x_rI_{m_r},0_k,-x_r I_{m_r},\ldots,-x_1
I_{m_1}\}|{x_i\in \bR}\},
$$
where $k=2m$ or $2m+1$ as appropriate.
We may denote an element of $\frak a$ by $x=(x_1,\ldots,x_r).$
We further note that $C_i\colon (x_1,\dots,x_i,\dots,x_r) \mapsto
(x_1,\dots,-x_i,\dots,x_r).$
Thus, if $C=C_B,$ we have $\frak a_C=\{(x_1,\ldots,x_r)|x_i=0,\forall i\in
B\}.$
Consequently, $\frak a_C=\{0\}$ if and only if $B=\{1,\ldots,r\}.$\qed
\enddemo

We now give more explicit criteria for $i_{G,M}(\s)$ to have elliptic
components.
In order to do this, we need to consider three cases.  The first two
are handled in the next result, and are similar to
the results of\Lspace \Lcitemark 15\Rcitemark \Rspace{}.  We use the same type
of arguments found there.

\proclaim{Theorem 3.4}
\roster
\item"(a)"If $\bG=GO_{2n+1}$ or $GU_{2n+1},$ then the following hold:

\itemitem{(i)}Conditions (a)--(c) of Proposition 3.3 hold if and only if
$R\simeq \Z_2^r.$

\itemitem{(ii)}If $\pi\in\calE_t(G),$ then either $\pi\in \calE_e (G),$ or
there is a proper Levi subgroup $\bM'$ of $\bG$ and some $\tau'\in\calE_e(M')$
with $\pi=i_{G,M'}(\tau).$

\item"(b)"If $\bG=GU_{2n}$ then the following hold:

\itemitem{(i)}Conditions (a)--(c) of Proposition 3.3 hold if and only if
$d=r-1$ and $d_{\omega_{E/F}} > 0$ is even, or $d=r$ and $d_{\omega_{E/F}}=0.$

\itemitem{(ii)}If $d < r-1,$ or $d=r-1$ and $d_{\omega_{E/F}}=1,$ then each
irreducible subrepresentation $\pi$ of $i_{G,M}(\sigma)$ is of the form
$\pi=i_{G,M'}(\tau),$ for some proper Levi subgroup $\bM'$ of $\bG,$ and some
$\tau\in\calE_t (M').$
We can choose $\bM'$ and $\tau$ with $\tau\in \calE_e (M')$
if and only if $d_{\omega_{E/F}}$ is even or $d_{\omega_{E/F}}=1.$

\itemitem{(iii)}Suppose $d=r-1,\ d_{\omega_{E/F}} \geq 3$ is odd, and $\pi$ is
an irreducible subrepresentation of $i_{G,M}(\sigma).$
Then $\pi$ cannot be irreducibly induced from a tempered representation of a
proper Levi subgroup.
\endroster
\endproclaim

\demo{Proof} (a) From Lemma 2.5, we see that if $\bG=GO_{2n+1}$ or $GU_{2n+1},$
then $C_1\ldots C_r\in R$ if and only if $R=\langle C_i|1\leq i\leq
r\rangle\simeq \Z_2^r.$
Thus, (i) holds.
Further, Lemma 2.5 implies that for any $\sigma\in \calE_2(M),$ there is a
subset $B_\sigma$ of $\{1,2,\ldots,r\}$ with $R=\langle C_i|i\in
B_\sigma\rangle.$
Let $C=C_{B_\sigma}.$
Then
$$\frak a_R=\{(x_1,\ldots,x_r)|x_i=0,\forall i\in B_\sigma\}=\frak a_C.$$
Thus, Theorem 1.5 implies (ii).

(b) Let $d_2=d_{\omega_{E/F}}.$
By the remark following Theorem 2.6, we know that $R\simeq \Z_2^{d_1+d_2-1}$ if
$d_2 > 0,$ and $R\simeq \Z_2^{d_1}$ if $d_2=0.$
So, if $d_2=0,$ we have $C_1\ldots C_r\in R$ if and only if $C_i\in R$ for each
$1\leq i\leq r,$
which is equivalent to $R\simeq \Z_2^r.$
On the other hand, if $d_2 > 0,$ then the longest element of $R$ consists of
$d_1+2[{d_2\over 2}]$ block sign changes.
Hence, if $d_2$ is odd, it is impossible for $C_1\ldots C_r\in R.$
On the other hand, if $d_2$ is even, then we see that $C_1 \ldots C_r\in R$ if
and only if $d_1+d_2=r,$ which is equivalent to $R\simeq \Z_2^{r-1}.$
This proves (i).

We note that if $d < r-1,$ then there is some $i\in\{1,\ldots,r\}$ for which
$C_i C_B \not\in R$ for all $B\subseteq \{1,\ldots,i-1,i+1,\ldots,r\}.$
Thus, for each $w\in R,$
the element $ E_i=(0,\ldots,1,\ldots,0)$
is in $\frak a_w.$ (Here the 1 is in the $i$\snug-th coordinate.)
Thus, $E_i \in \frak a_R,$ which implies $\frak a_R\not= \{0\}.$
Therefore, by Theorem 1.5, every irreducible subrepresentation $\pi$ of
$i_{G,M}(\sigma)$ is of the form $\pi=i_{G,M'}(\tau)$ for some $\M'\subsetneq
\bG$ and $\tau\in \calE_t(M').$

Suppose now that $d_2=d_{\omega_{E/F}}$ is even.  Let
$I_2(\s)=I_{\o_{E/F}}(\s).$
Then
$$C_0=\prod\limits_{i\in I_2(\sigma)} C_i$$
is in $R.$
So if $B=I_1 (\sigma)\cup I_2(\s),$ we have $C_B \in R,$ and $\frak
a_{C_B}=\frak a_R.$
Therefore, if $d < r-1$ and $d_2$ is even, then we can choose $(\bM',\tau),$
with $\tau \in \calE_e (M').$

Suppose that $d_2=1.$
Then $\o_{E/F}\not\in X(\rho).$
If $\o_1=\o_{E/F},$
then $C_1C_B\not\in R,$ for all $B\subset\{2,\dots,r\}.$
Let $C=C_{I_1(\s)}.$ Then $\frak a_R=\frak a_C,$
and therefore we can choose $(\M',\tau)$ and $\tau\in\Cal E_e(M'),$
with $\pi=i_{G,M'}(\tau).$

Now suppose $d_2\geq 3$ is odd.  Without loss of generality, we assume
that $I_1(\s)=\{k+1,k+2,\dots,k+d_1\},$ and
$I_2(\sigma)=\{k+d_1+1,\ldots,r-1,r\},$ with $k=r-(d_1+d_2).$
Then $\frak a_R=\{(x_1 \ldots x_r)|x_{k+1}=x_{k+2}=\ldots=x_r=0\}.$
Since $C_{k+1} \ldots C_r\not\in R,$ we see that $\frak a_R \not= \frak a_w$
for any $w\in R.$
Therefore, $\pi$ can not be of the form $i_{G,M'}(\tau)$ for some $\bM'
\subsetneq \bG$ and $\tau\in \calE_e (M').$
\qed
\enddemo

\definition{Definition 3.5}
Let $\G=GSp_{2n}$ or $GO_{2n}.$
Suppose $\s\simeq\s_1\otimes\dots\otimes\s_r\otimes\rho\in\Cal E_2(M).$
We set
$$O(G,\s)=\left\{\chi\in\left(F^\times/F^{\times^2}\right)^\wedge\setminus
X(\rho)\ \big |\ d_{\chi}\text{ is odd }\right \}.$$
We also set
$$O_1(G,\s)=\left\{\chi\in O(G,\s)\ \big|\ d_\chi=1,\text{ and }\chi\not\in S,\
\forall S\in\L(\s)\right\}.$$
Finally, we let
$I_0(\s)=I_1(\s)\cup\bigcup\limits_{\chi\not\in X(\rho)}I_{\chi}(\s).$
\enddefinition

\proclaim{Theorem 3.6} Suppose $\bG=GSp_{2n}$ or $GO_{2n}.$
Let $M\simeq GL_{m_1}\times \ldots\times GL_{m_r}\times G(m).$
Assume that
$\sigma\simeq\s_1\otimes\dots\otimes\s_r\otimes\rho\in\calE_2 (M)$
and $R$ is the $R$\snug-group of $\sigma.$

\roster
\item"(a)"$i_{G,M}(\sigma)$ has elliptic constituents if and only if

\itemitem{(i)}$O(G,\sigma)$ can be partitioned into minimally
$\rho$\snug-trivial subsets, and

\itemitem{(ii)}$I_0(\s)=\{1,\ldots,r\}.$

\item"(b)" Suppose $\pi\in \calE_t (G)$ is not elliptic, and $\pi\subseteq
i_{G,M}(\sigma).$  Then we can find a proper Levi subgroup $\bM'$ of $ \bG$
and a $\tau\in \calE_t (M')$ with $\pi=i_{G,M'}(\tau)$
if and only if
$I_0(\s) \subsetneq \{1,\ldots,r\},$
or $O_1(G,\s)$ is nonempty.
We can choose $(\bM',\tau)$ with $\tau\in \calE_e (M')$ if and only if
$O(G,\sigma)\backslash O_1(G,\sigma)$ can be partitioned into minimally
$\rho$\snug-trivial subsets.
\endroster
\endproclaim

\demo{Proof}
(a) If $C_1\ldots C_r \in R,$ then $I_0 (\sigma)=\{1,\ldots,r\}.$
So condition (ii) is trivial.
Suppose this is the case.
For each $\chi\in O(G,\sigma),$ we choose an $i_\chi \in \{1,\ldots,r\},$ with
$\omega_{i_\chi}=\chi.$
We then let
$$
C_\chi=\prod_{i\in I_\chi (\sigma)\backslash \{i_\chi\}} C_i.
$$
Since $d_\chi$ is odd, $C_\chi\in R.$
For each $\chi\not\in X(\rho)$ with $d_\chi$ even, we let
$$C_\chi=\prod\limits_{i\in I_\chi (\sigma)} C_i.$$
Finally let
$$
C_0=\left(\prod_{i\in I_1 (\sigma)} C_i\right) \prod_{\chi\not\in X(\rho)}
C_\chi.
$$
Then $C_0\in R,$ and
$$C_0=\prod\limits_{i\not\in \{i_\chi|\chi\in O(G,\sigma)\}}C_i.$$
Thus, $C_1 \ldots C_r\in R$ if and only if
$C'=\prod\limits_{O(G,\sigma)}C_{i_\chi} \in R.$
 From the proof of Theorem 2.6, we know that $C'\in R$ if and only if
$\prod\limits_{O(G,\sigma)} \omega_{i_\chi}=\prod\limits_{O(G,\sigma)}\chi\in
X(\rho).$
This shows that (i) and (ii) are equivalent to $C_1\ldots C_r\in R.$

(b) Note that
$$\frak a_R=\bigcap\limits_{w\in R}\frak a_w \subseteq \left\{(x_1\ldots x_r)\
\big|\ x_i=0,\ \forall i\in I_0(\sigma)\right\}.$$
So in particular, if $I_0 (\sigma)\subsetneq \{1,\ldots,r\},$ then $\frak a_R
\supsetneq\{0\}.$
Thus, in this case,
we can find $\bM'$ and $\tau\in \Cal E_t(M')$ as desired.
Now suppose that $I_0 (\sigma)=\{1,\ldots,r\}.$
If $\chi\not\in X(\rho)$ and $\omega_i=\omega_j=\chi$ for some $i\not= j,$ then
$C=C_i C_j \in R$ and $\frak a_C=\{(x_1\ldots x_r)|x_i=x_j=0\}\subseteq \frak
a_R.$
Therefore, if $d_\chi\not= 1$ for all $\chi\not\in X(\rho)$
then we see that $\frak a_R=\{0\},$ and so $i_{G,M'}(\tau)=\pi$ is impossible.
So now suppose $\omega_1=\chi$ and $\omega_i\not= \chi$ for all $i > 1.$
Then there is a subset $B\subset \{2,3,\ldots,r\}$ with $C_1 C_B \in R$ if and
only if $\o_B=\chi,$ which is equivalent to $\{\chi\} \cup \{\omega_j|j\in
B'\}\in \Lambda (\sigma)$ for some $B'\subseteq \{2,\dots,r\}.$
Thus, if no such $B'$ exists then $\{(1,0,\ldots,0)\}\in \frak a_R,$ and we can
choose $\M'$ and $\tau\in \calE_t (M')$ with $\pi=i_{G,M'}(\tau).$
On the other hand, suppose that, for each $\chi$ with $d_\chi=1,$ there is an
$S\in\Lambda(\sigma)$ with $\chi\in S.$
Then, for each $1\leq i\leq r,$ there is a subset $B$ of
$\{1,\ldots,i-1,i+1,\ldots,r\}$ with $C_i C_B\in R.$
Thus, for each $i,$ $\{(x_1\ldots x_r)|x_i=0\}\subseteq \frak a_R,$ which
implies $\frak a_R=\{0\}.$
Therefore, finding $\M'$ and $\tau$ is impossible.

Finally, we assume that $\frak a_R=\{(x_1,\ldots,x_k,0,\ldots,0)|x_i\in \bR\},$
with $k\geq 1.$
We need to determine when $C_{k+1}\ldots C_r\in R.$
Let $\chi\in O(G,\sigma).$
If $\chi\in O_1(G,\s),$
and $\o_{i_\chi}=\chi,$ then $i_\chi \leq k.$
For all other $\chi\in O(G,\sigma),$ if $i\in I_\chi(\s),$ then
we have $i > k.$
Suppose $\chi\in O(G,\sigma)\backslash O_1(G,\sigma).$
As in the proof of part (a),  we fix $i_\chi \in I_\chi (\sigma)$ and let
$$C_\chi=\prod\limits_{i\in I_\chi (\sigma)\backslash \{ i_\chi\}} C_i.$$
If $\chi\not\in O(G,\sigma)$ let
$$C_\chi=\prod\limits_{i\in I_\chi(\sigma)}C_i.$$
Set
$$C_0=\prod\limits_{i\in I_1(\s)}C_i \prod\limits_{\chi\not\in
X(\rho)}C_\chi.$$
Then $C_0\in R.$
Since
$$C_{k+1}\ldots C_r=C_0\cdot \prod\limits_{\chi\in O(G,\sigma)\backslash
O_1(G,\sigma)} C_{i_\chi},$$
we see that $C_{k+1}\ldots C_r \in R$ if and only if
$$\prod\limits_{O(G,\sigma)\backslash O_1(G,\sigma)}\chi\in X(\rho),$$
which says that $O(G,\s)\setminus O_1(G,\s)$
can be partitioned into minimally $\rho$\snug-trivial subsets.\qed
\enddemo

\example{Example}
We note that in all previous cases where the elliptic tempered
representations were determined, the following phenomenon always held:
Given a parabolic subgroup $\PPP=\M\N$ of $\G,$  there is a fixed finite group
$R_M,$ so that, for any $\s\in\Cal E_2(M),$
the representation $i_{G,M}(\s)$ has elliptic components if
and only if $R(\s)\simeq R_M.$  In fact, if $\M$ is a Levi subgroup that
admits elliptic tempered representations, then $i_{G,M}(\s)$ has
elliptic constituents if and only if $|R(\s)|=|R_M|.$  With the similitude
groups, we see that such $R_M$ do not exist.  We look at
the following example.  Suppose that $G=GSp_{6n}$ or $GO_{6n}$
for some $n,$ and $\M\simeq GL_n\times GL_n\times GL_n\times GL_1.$
Let $\s\simeq\s_1\otimes\s_2\otimes\s_3\otimes\rho.$
In this case $\rho$ is a one dimensional character,
so $\rho\otimes\o=\rho$ if and only if
$\o=1.$  We suppose that
the $\s_i$ are pairwise inequivalent discrete series representations
of $GL_n(F),$ and that
$\s_i\simeq\wt\s_i\simeq\s_i^\e$ for each $i.$
Here $\wt\s_i$ is the smooth contragredient of $\s_i$
\Lcitemark 3\LIcitemark{},\Sec 7\RIcitemark \Rcitemark \Rspace{}.
If $\Cal X_{n,0,G}(\s_i\otimes\rho)$ holds for each $i,$ then
$R\simeq\Z_2\times\Z_2\times\Z_2,$ and $i_{G,M}(\s)$ has eight elliptic
constituents.
On the other hand, if $\Cal X_{n,0,G}(\s_1\otimes\rho)$ holds, while
$\o_2=\o_3\neq 1,$
then $R=<C_1,C_2C_3>.$  Since $C_1C_2C_3\in R,$  we see that
$i_{G,M}(\s)$ has four elliptic constituents.  Finally, if the $\o_i$
are all non-trivial, and $\o_3=\o_1\o_2,$ then $R=<C_1C_2C_3>,$
and $i_{G,M}(\s)$ has two elliptic constituents.
\qed\endexample

Suppose now that $\bG$ is any of the similitude groups that we have been
studying.
Further suppose $\bM$ is a Levi subgroup of $\bG,$ and $\sigma\in \calE_2 (M)$
is such that all the constituents of $i_{G,M}(\sigma)$ are elliptic.
Let $\hat{R}$ be the collection of irreducible representations of $R.$
Since $R$ is abelian, $\dim \k=1$ for all $\k\in \hat{R}.$
We let $\pi_\kappa$ be the irreducible subrepresentation of $i_{G,M}(\sigma)$
canonically parameterized by $\kappa$\Lspace \Lcitemark 19\Rcitemark \Rspace{},
and we
denote by $\Theta^e_\k$
its character on the regular elliptic set of $G.$
Then we see that
$$
\sum_{\kappa\in \hat{R}} \Theta_\kappa^e=0.
$$
We now make this linear dependence precise by generalizing an argument of
Herb\Lspace \Lcitemark 15\Rcitemark \Rspace{}.
We need some preliminary results.

\proclaim{Lemma 3.7} Suppose $i_{G,M}(\sigma)$ has elliptic constituents, and
$R\simeq \Z_2^d.$
For $1\leq i\leq r,$ let $\bM_i$ be the standard Levi component of $\bG,$
with split component $\A_i,$
such that $M_i\simeq GL_{m_i}\times G(n-m_i),$
and $W(\bG,\bA_i)=\langle C_i\rangle. $
Then the $R$\snug-group $R_i$ attached to $i_{M_i,M}(\sigma)$ is of order
$\Z_2^{d-1}.$
\endproclaim

\demo{Proof} Since $i_{G,M}(\sigma)$ has elliptic constituents, $C_1\ldots
C_r\in R,$ and thus $\Delta'=\emptyset.$
Therefore, $\bM_i$ satisfies Arthur's compatibility condition, which implies
$R_i=R\cap W(\bM_i,\bA).$
Suppose $R=\langle s _1,\ldots,s_d\rangle.$
We assume that $s_1=C_i C_B$ for some $B\subseteq
\{1,\ldots,i-1,i+1,\ldots,r\}.$
We claim that it is possible to choose $\{s_1,\ldots,s_d\}$ so that $s_1$ is
the only generator of this form.
For $j>1,$ we suppose that $s_j=C_{B_j}.$
Set $s_j'=s_j$ if $i\not\in B_j,$ and
$s_j'=s_1s_j$
otherwise.
Then  $\{s_1,s_2',\ldots,s'_d\}$ is a set of generators of
the desired form.  Therefore, $R_i=R\cap W(\bM_i, \bA)=\langle
s'_2,\ldots,s'_d\rangle \simeq \Z_2^{d-1}.$
\qed
\enddemo

\proclaim{Lemma 3.8}
Suppose that $M\simeq GL_{m_1}\times\dots\times GL_{m_r}\times G(m).$
Let $\s\in\Cal E_2(M),$ and suppose that $i_{G,M}(\s)$
has elliptic constituents.  We further suppose that
$R=R(\s)\simeq\Z_2^d.$  Then there is a set of
generators $\O=\{s_1,\dots,s_d\}$ for $R$ so that,
for each $i,$
$\O\setminus\{s_i\}$ generates $R_{j_i},$ for some $j_i.$
\endproclaim

\demo{Proof}
By Lemma 3.7 we can choose
a set of generators $\{s_{11},\dots,s_{1d}\}$ for $R$
with
$R_{j_1}=<s_{12},\dots,s_{1d}>$ for some $j_1.$
Now suppose $1\leq k<d,$ and we have chosen
a generating set $\O_k=\{s_{k1},\dots,s_{kd}\}$ for $R$
with the property that, for each $1\leq i\leq k$
the set
$\O_k\setminus\{s_{ki}\}$
generates $R_{j_i}$ for some $j_i.$
Suppose $s_{k(k+1)}=C_B,$ for some $B\subset\{1,\dots,r\}.$
Since
$$s_{k(k+1)}\in\bigcap_{i=1}^k R_{j_i},$$
we know that $j_i\not\in B$ for each $1\leq i\leq k.$
Fix $j_{k+1}\in B.$
Suppose that $i\neq k+1$ and
$s_{ki}=C_{B_i}.$
For $i\neq k+1,$ we let
$$s_{(k+1)i}=\cases s_{ki}&\text{ if }j_{k+1}\not\in B_i\\
s_{k(k+1)}s_{ki}&\text{ if }j_{k+1}\in B_i,\endcases$$
and take $s_{(k+1)(k+1)}=s_{k(k+1)}.$
Set $\O_{k+1}=\{s_{(k+1)1},\dots,s_{(k+1)d}\}.$
By our choice of
$\O_{k},\ \O_{k+1},$ and by Lemma 3.7,
we see that for $1\leq i\leq k,$ the group
$R_{j_i}$ is generated by $\O_{k+1}\setminus\{s_{(k+1)i}\}.$
Moreover, we have chosen the $s_{(k+1)i}=C_{D_i}$ so that
$j_{k+1}\not\in D_i,$ for $i\neq k+1.$  Thus,
$R_{j_{k+1}}$ is generated by $\O_{k+1}\setminus\{s_{(k+1)(k+1)}\}.$
Consequently, by induction, we can choose $\O=\O_d$ with the desired
property.\qed
\enddemo

\proclaim{Theorem 3.9} Let $G=GSp_{2n},\ GO_n,$ or $GU_n.$
Suppose $\M$ is a Levi subgroup of $\G,$ with $M\simeq GL_{m_1} \times \ldots
\times GL_{m_r}\times G(m).$
Further suppose that $\sigma\in \Cal E_2 (M)$ satisfies conditions (a)--(c) of
Proposition 3.3.
Let $R=R(\sigma).$
Suppose $\kappa\in \hat{R}$ and that $\pi_\kappa$ is the irreducible
subrepresentation of $i_{G,M}(\sigma)$ parameterized by $\kappa.$ Denote its
character on $G^e$ by $\Theta^e_\kappa.$
We set $w_0=C_1\dots C_r,$ and let
$\e(\k)=\k(w_0).$
Then $\Theta_\kappa^e=\e(\k)\Theta_1^e,$
where $1$ is the trivial character.
\endproclaim

\demo{Proof}
By Lemma 3.8, we can choose a set of generators
$\O=\{s_1,\dots,s_d\}$ for $R$ with the
property that $\O\setminus\{s_i\}$ generates $R_{j_i}$
for some $j_i.$  We fix a choice of such an $\O.$
Let $\k\in\hat R.$  Denote by $s(\k)$
the number of $s_i\in\O$ for which $\k(s_i)=-1.$
If $s(\k)=0,$ then $\k=1,$ and the claim is trivially true.
Assume that the statement of the theorem is true
whenever $s(\k)\leq s.$  Suppose that $s(\k)=s+1.$  Without
loss of generality, we suppose that $\k(s_1)=-1.$
We are assuming that $<s_2,\dots,s_d>= R_j$ for some $j.$
Let $\k_j=\k|_{{R_j}}.$
If $\xi\in\hat R_j,$ then
$$\hat R(\xi)=\{\chi\in\hat R\ |\ \chi|_{R_j}=\xi\}=\{\xi^+,\xi^-\},$$
where $\xi^{\pm}$ is  determined by $\xi^{\pm}(s_1)=\pm1.$
Note that $\k=\k_j^-,$ and consequently
$s(\k_j^+)=s(\k)-1=s.$  Therefore, by our induction
hypothesis
$\T_{\k_j^+}^e=\e(\k_j^+)\T_1^e.$
Let $\tau_{\k_j}$ be the irreducible subrepresentation
of $i_{M_j,M}(\s)$ parameterized by $\k_j.$
Then, by Theorem 1.4, we have
$i_{G,M_j}(\tau_{\k_j})=\pi_{\k_j^+}\oplus\pi_{\k_j^-},$
and thus,
$\T_\k^e=-\T_{\k_j^+}^e=-\e(\k_j^+)\T_1^e.$
It is enough to show that $\e(\k)=-\e(\k_j^+).$
To see this, write $w_0=s_1w,$ with $w\in R_j.$
Now we have
$$\e(\k)=\k(w_0)=\k(s_1)\k(w)=-\k_j^+(s_1)\k_j^+(w)=-\k_j^+(w_0)=-\e(\k_j^+).$$
Therefore, $\T_\k^e=\e(\k)\T_1^e,$ and the theorem follows
by induction.\qed
\vfil\eject
\enddemo

\Refs

\message{REFERENCE LIST}

\bgroup\Resetstrings%
\def\Ecnt{0}\def\acnt{0}%
\def\Ftest{ }\def\Fstr{1}%
\def\Atest{ }\def\Astr{J\Initper  Arthur}%
\def\Ttest{ }\def\Tstr{Intertwining operators and residues I. Weighted
characters}%
\def\Jtest{ }\def\Jstr{J. Funct. Anal.}%
\def\Vtest{ }\def\Vstr{84}%
\def\Ptest{ }\def\Pstr{19--84}%
\def\Dtest{ }\def\Dstr{1989}%
\Refformat\egroup%

\bgroup\Resetstrings%
\def\Ecnt{0}\def\acnt{0}%
\def\Ftest{ }\def\Fstr{2}%
\def\Atest{ }\def\Astr{J\Initper  Arthur}%
\def\Ttest{ }\def\Tstr{On elliptic tempered characters}%
\def\Jtest{ }\def\Jstr{Acta Math.}%
\def\Vtest{ }\def\Vstr{171}%
\def\Dtest{ }\def\Dstr{1993}%
\def\Ptest{ }\def\Pstr{73--138}%
\def\Astr{\Underlinemark}%
\Refformat\egroup%

\bgroup\Resetstrings%
\def\Ecnt{0}\def\acnt{0}%
\def\Ftest{ }\def\Fstr{3}%
\def\Atest{ }\def\Astr{I\Initper \Initgap N\Initper  Bernstein%
  \Aand A\Initper \Initgap V\Initper  Zelevinsky}%
\def\Ttest{ }\def\Tstr{Representations of the group $GL(n,F)$ where $F$ is a
local non-archimedean local field}%
\def\Jtest{ }\def\Jstr{Russian Math. Surveys}%
\def\Vtest{ }\def\Vstr{33}%
\def\Ptest{ }\def\Pstr{1--68}%
\def\Dtest{ }\def\Dstr{1976}%
\Refformat\egroup%

\bgroup\Resetstrings%
\def\Ecnt{0}\def\acnt{0}%
\def\Ftest{ }\def\Fstr{4}%
\def\Atest{ }\def\Astr{I\Initper \Initgap N\Initper  Bernstein%
  \Aand A\Initper \Initgap V\Initper  Zelevinsky}%
\def\Ttest{ }\def\Tstr{Induced representations of reductive $p$\snug-adic
groups. I}%
\def\Jtest{ }\def\Jstr{Ann. Sci. \'Ecole Norm. Sup. (4)}%
\def\Vtest{ }\def\Vstr{10}%
\def\Ptest{ }\def\Pstr{441--472}%
\def\Dtest{ }\def\Dstr{1977}%
\def\Astr{\Underlinemark}%
\Refformat\egroup%

\bgroup\Resetstrings%
\def\Ecnt{0}\def\acnt{0}%
\def\Ftest{ }\def\Fstr{5}%
\def\Atest{ }\def\Astr{N\Initper  Bourbaki}%
\def\Ttest{ }\def\Tstr{Groupes et Alg\`ebres de Lie, Chapitres 4,5, et 6}%
\def\Itest{ }\def\Istr{Hermann}%
\def\Ctest{ }\def\Cstr{Paris}%
\def\Dtest{ }\def\Dstr{1968}%
\Refformat\egroup%

\bgroup\Resetstrings%
\def\Ecnt{0}\def\acnt{0}%
\def\Ftest{ }\def\Fstr{6}%
\def\Atest{ }\def\Astr{D\Initper  Goldberg}%
\def\Ttest{ }\def\Tstr{$R$\snug-groups and elliptic representations for unitary
groups}%
\def\Jtest{ }\def\Jstr{Proc. Amer. Math. Soc.}%
\def\Otest{ }\def\Ostr{to appear}%
\Refformat\egroup%

\bgroup\Resetstrings%
\def\Ecnt{0}\def\acnt{0}%
\def\Ftest{ }\def\Fstr{7}%
\def\Atest{ }\def\Astr{D\Initper  Goldberg}%
\def\Ttest{ }\def\Tstr{Reducibility of induced representations for $Sp(2n)$ and
$SO(n)$}%
\def\Jtest{ }\def\Jstr{Amer. J. Math.}%
\def\Otest{ }\def\Ostr{to appear}%
\def\Astr{\Underlinemark}%
\Refformat\egroup%

\bgroup\Resetstrings%
\def\Ecnt{0}\def\acnt{0}%
\def\Ftest{ }\def\Fstr{8}%
\def\Atest{ }\def\Astr{D\Initper  Goldberg}%
\def\Ttest{ }\def\Tstr{Reducibility for non-connected $p$\snug-adic groups with
$G^\circ$ of prime index}%
\def\Jtest{ }\def\Jstr{Canad. J. Math.}%
\def\Otest{ }\def\Ostr{to appear}%
\def\Astr{\Underlinemark}%
\Refformat\egroup%

\bgroup\Resetstrings%
\def\Ecnt{0}\def\acnt{0}%
\def\Ftest{ }\def\Fstr{9}%
\def\Atest{ }\def\Astr{D\Initper  Goldberg}%
\def\Ttest{ }\def\Tstr{Some results on reducibility for unitary groups and
local Asai $L$\snug-functions}%
\def\Jtest{ }\def\Jstr{J. Reine Angew. Math.	}%
\def\Vtest{ }\def\Vstr{448}%
\def\Ptest{ }\def\Pstr{65--95}%
\def\Dtest{ }\def\Dstr{1994}%
\def\Astr{\Underlinemark}%
\Refformat\egroup%

\bgroup\Resetstrings%
\def\Ecnt{0}\def\acnt{0}%
\def\Ftest{ }\def\Fstr{10}%
\def\Atest{ }\def\Astr{D\Initper  Goldberg}%
\def\Ttest{ }\def\Tstr{$R$\snug-groups and elliptic representations for
$SL_n$}%
\def\Jtest{ }\def\Jstr{Pacific J. Math.}%
\def\Vtest{ }\def\Vstr{165}%
\def\Dtest{ }\def\Dstr{1994}%
\def\Ptest{ }\def\Pstr{77--92}%
\def\Astr{\Underlinemark}%
\Refformat\egroup%

\bgroup\Resetstrings%
\def\Ecnt{0}\def\acnt{0}%
\def\Ftest{ }\def\Fstr{11}%
\def\Atest{ }\def\Astr{D\Initper  Goldberg%
  \Aand F\Initper  Shahidi}%
\def\Ttest{ }\def\Tstr{Automorphic $L$\snug-functions, intertwining operators,
and the irreducible tempered representations of $p$\snug-adic groups}%
\def\Otest{ }\def\Ostr{to appear}%
\Refformat\egroup%

\bgroup\Resetstrings%
\def\Ecnt{0}\def\acnt{0}%
\def\Ftest{ }\def\Fstr{12}%
\def\Atest{ }\def\Astr{D\Initper  Goldberg%
  \Aand F\Initper  Shahidi}%
\def\Ttest{ }\def\Tstr{On the tempered spectrum of quasi-split classical
groups}%
\def\Otest{ }\def\Ostr{preprint}%
\def\Astr{\Underlinemark}%
\Refformat\egroup%

\bgroup\Resetstrings%
\def\Ecnt{0}\def\acnt{0}%
\def\Ftest{ }\def\Fstr{13}%
\def\Atest{ }\def\Astr{Harish-Chandra}%
\def\Ttest{ }\def\Tstr{Harmonic Analysis on Reductive $p$\snug-adic Groups}%
\def\Rtest{ }\def\Rstr{Notes by G. van Dijk}%
\def\Itest{ }\def\Istr{Springer-Verlag}%
\def\Ctest{ }\def\Cstr{New York--Heidelberg--Berlin}%
\def\Stest{ }\def\Sstr{Lecture Notes in Mathematics}%
\def\Ntest{ }\def\Nstr{162}%
\def\Dtest{ }\def\Dstr{1970}%
\Refformat\egroup%

\bgroup\Resetstrings%
\def\Ecnt{0}\def\acnt{0}%
\def\Ftest{ }\def\Fstr{14}%
\def\Atest{ }\def\Astr{Harish-Chandra}%
\def\Ttest{ }\def\Tstr{Harmonic analysis on reductive $p$\snug-adic groups}%
\def\Dtest{ }\def\Dstr{1973}%
\def\Jtest{ }\def\Jstr{Proc. Sympos. Pure Math.}%
\def\Itest{ }\def\Istr{AMS}%
\def\Ctest{ }\def\Cstr{Providence, RI}%
\def\Vtest{ }\def\Vstr{26}%
\def\Ptest{ }\def\Pstr{167--192}%
\def\Astr{\Underlinemark}%
\Refformat\egroup%

\bgroup\Resetstrings%
\def\Ecnt{0}\def\acnt{0}%
\def\Ftest{ }\def\Fstr{15}%
\def\Atest{ }\def\Astr{R\Initper \Initgap A\Initper  Herb}%
\def\Ttest{ }\def\Tstr{Elliptic representations for $Sp(2n)$ and $SO(n)$}%
\def\Jtest{ }\def\Jstr{Pacific J. Math.}%
\def\Vtest{ }\def\Vstr{161}%
\def\Dtest{ }\def\Dstr{1993}%
\def\Ptest{ }\def\Pstr{347--358}%
\Refformat\egroup%

\bgroup\Resetstrings%
\def\Ecnt{0}\def\acnt{0}%
\def\Ftest{ }\def\Fstr{16}%
\def\Atest{ }\def\Astr{N\Initper  Jacobson}%
\def\Ttest{ }\def\Tstr{A note on hermitian forms}%
\def\Jtest{ }\def\Jstr{Bull. Amer. Math. Soc.}%
\def\Vtest{ }\def\Vstr{46}%
\def\Dtest{ }\def\Dstr{1940}%
\def\Ptest{ }\def\Pstr{264--268}%
\Refformat\egroup%

\bgroup\Resetstrings%
\def\Ecnt{0}\def\acnt{0}%
\def\Ftest{ }\def\Fstr{17}%
\def\Atest{ }\def\Astr{H\Initper  Jacquet}%
\def\Ttest{ }\def\Tstr{Generic representations}%
\def\Btest{ }\def\Bstr{Non Commutatuve Harmonic Analysis}%
\def\Itest{ }\def\Istr{Springer-Verlag}%
\def\Ctest{ }\def\Cstr{New York--Heidelberg--Berlin}%
\def\Stest{ }\def\Sstr{Lecture Notes in Mathematics}%
\def\Ntest{ }\def\Nstr{587}%
\def\Dtest{ }\def\Dstr{1977}%
\def\Ptest{ }\def\Pstr{91--101}%
\Refformat\egroup%

\bgroup\Resetstrings%
\def\Ecnt{0}\def\acnt{0}%
\def\Ftest{ }\def\Fstr{18}%
\def\Atest{ }\def\Astr{C\Initper \Initgap D\Initper  Keys}%
\def\Ttest{ }\def\Tstr{On the decomposition of reducible principal series
representations of $p$\snug-adic Chevalley groups}%
\def\Dtest{ }\def\Dstr{1982}%
\def\Jtest{ }\def\Jstr{Pacific J. Math.}%
\def\Vtest{ }\def\Vstr{101}%
\def\Ptest{ }\def\Pstr{351--388}%
\Refformat\egroup%

\bgroup\Resetstrings%
\def\Ecnt{0}\def\acnt{0}%
\def\Ftest{ }\def\Fstr{19}%
\def\Atest{ }\def\Astr{C\Initper \Initgap D\Initper  Keys}%
\def\Ttest{ }\def\Tstr{L-indistinguishability and R-groups for quasi split
groups: unitary groups in even dimension}%
\def\Jtest{ }\def\Jstr{Ann. Sci. \'Ecole Norm. Sup. (4)}%
\def\Dtest{ }\def\Dstr{1987}%
\def\Vtest{ }\def\Vstr{20}%
\def\Ptest{ }\def\Pstr{31--64}%
\def\Astr{\Underlinemark}%
\Refformat\egroup%

\bgroup\Resetstrings%
\def\Ecnt{0}\def\acnt{0}%
\def\Ftest{ }\def\Fstr{20}%
\def\Atest{ }\def\Astr{A\Initper \Initgap W\Initper  Knapp%
  \Aand E\Initper \Initgap M\Initper  Stein}%
\def\Ttest{ }\def\Tstr{Irreducibility theorems for the principal series}%
\def\Btest{ }\def\Bstr{Conference on Harmonic Analysis}%
\def\Itest{ }\def\Istr{Springer-Verlag}%
\def\Ctest{ }\def\Cstr{New York--Heidelberg--Berlin}%
\def\Stest{ }\def\Sstr{Lecture Notes in Mathematics}%
\def\Dtest{ }\def\Dstr{1972}%
\def\Ntest{ }\def\Nstr{266}%
\def\Ptest{ }\def\Pstr{197--214}%
\Refformat\egroup%

\bgroup\Resetstrings%
\def\Ecnt{0}\def\acnt{0}%
\def\Ftest{ }\def\Fstr{21}%
\def\Atest{ }\def\Astr{A\Initper \Initgap W\Initper  Knapp%
  \Aand G\Initper  Zuckerman}%
\def\Ttest{ }\def\Tstr{Classification of irreducible tempered representations
of semisimple Lie groups}%
\def\Jtest{ }\def\Jstr{Proc. Nat. Acad. Sci. U.S.A.}%
\def\Vtest{ }\def\Vstr{73}%
\def\Ntest{ }\def\Nstr{7}%
\def\Ptest{ }\def\Pstr{2178--2180}%
\def\Dtest{ }\def\Dstr{1976}%
\Refformat\egroup%

\bgroup\Resetstrings%
\def\Ecnt{0}\def\acnt{0}%
\def\Ftest{ }\def\Fstr{22}%
\def\Atest{ }\def\Astr{W\Initper  Landherr}%
\def\Ttest{ }\def\Tstr{\"Aquivalenze Hermitscher formen \"uber einem beliebigen
algebraischen zahlk\"orper}%
\def\Jtest{ }\def\Jstr{Abh. Math. Sem. Univ. Hamburg}%
\def\Vtest{ }\def\Vstr{11}%
\def\Dtest{ }\def\Dstr{1936}%
\def\Ptest{ }\def\Pstr{245--248}%
\Refformat\egroup%

\bgroup\Resetstrings%
\def\Ecnt{0}\def\acnt{0}%
\def\Ftest{ }\def\Fstr{23}%
\def\Atest{ }\def\Astr{G\Initper \Initgap I\Initper  Ol'{\v{s}}anski{\v{i}}}%
\def\Ttest{ }\def\Tstr{Intertwining operators and complementary series in the
class of representations induced from parabolic subgroups of the general linear
group over a locally compact division algebra}%
\def\Jtest{ }\def\Jstr{Math. USSR-Sb.}%
\def\Dtest{ }\def\Dstr{1974}%
\def\Vtest{ }\def\Vstr{22}%
\def\Ptest{ }\def\Pstr{217--254}%
\Refformat\egroup%

\bgroup\Resetstrings%
\def\Ecnt{0}\def\acnt{0}%
\def\Ftest{ }\def\Fstr{24}%
\def\Atest{ }\def\Astr{F\Initper  Shahidi}%
\def\Ttest{ }\def\Tstr{The notion of norm and the representation theory of
orthogonal groups}%
\def\Jtest{ }\def\Jstr{Invent. Math.}%
\def\Otest{ }\def\Ostr{To appear}%
\Refformat\egroup%

\bgroup\Resetstrings%
\def\Ecnt{0}\def\acnt{0}%
\def\Ftest{ }\def\Fstr{25}%
\def\Atest{ }\def\Astr{F\Initper  Shahidi}%
\def\Ttest{ }\def\Tstr{On certain $L$\snug-functions}%
\def\Jtest{ }\def\Jstr{Amer. J. Math.}%
\def\Vtest{ }\def\Vstr{103}%
\def\Ntest{ }\def\Nstr{2}%
\def\Ptest{ }\def\Pstr{297--355}%
\def\Dtest{ }\def\Dstr{1981}%
\def\Astr{\Underlinemark}%
\Refformat\egroup%

\bgroup\Resetstrings%
\def\Ecnt{0}\def\acnt{0}%
\def\Ftest{ }\def\Fstr{26}%
\def\Atest{ }\def\Astr{F\Initper  Shahidi}%
\def\Ttest{ }\def\Tstr{Some results on $L$\snug-indistinguishability for
$SL(r)$}%
\def\Jtest{ }\def\Jstr{Canad. J. Math.}%
\def\Dtest{ }\def\Dstr{1983}%
\def\Vtest{ }\def\Vstr{35}%
\def\Ptest{ }\def\Pstr{1075--1109}%
\def\Astr{\Underlinemark}%
\Refformat\egroup%

\bgroup\Resetstrings%
\def\Ecnt{0}\def\acnt{0}%
\def\Ftest{ }\def\Fstr{27}%
\def\Atest{ }\def\Astr{F\Initper  Shahidi}%
\def\Ttest{ }\def\Tstr{Fourier transforms of intertwining operators and
Placherel measures for  $GL(n)$}%
\def\Jtest{ }\def\Jstr{Amer. J. Math.}%
\def\Vtest{ }\def\Vstr{106}%
\def\Dtest{ }\def\Dstr{1984}%
\def\Ptest{ }\def\Pstr{67--111}%
\def\Astr{\Underlinemark}%
\Refformat\egroup%

\bgroup\Resetstrings%
\def\Ecnt{0}\def\acnt{0}%
\def\Ftest{ }\def\Fstr{28}%
\def\Atest{ }\def\Astr{F\Initper  Shahidi}%
\def\Ttest{ }\def\Tstr{A proof of Langlands conjecture for Plancherel measures;
complementary series for $p$\snug-adic groups}%
\def\Jtest{ }\def\Jstr{Ann. of Math. (2)}%
\def\Dtest{ }\def\Dstr{1990}%
\def\Vtest{ }\def\Vstr{132}%
\def\Ptest{ }\def\Pstr{273--330}%
\def\Astr{\Underlinemark}%
\Refformat\egroup%

\bgroup\Resetstrings%
\def\Ecnt{0}\def\acnt{0}%
\def\Ftest{ }\def\Fstr{29}%
\def\Atest{ }\def\Astr{F\Initper  Shahidi}%
\def\Ttest{ }\def\Tstr{Twisted endoscopy and reducibility of induced
representations for $p$\snug-adic groups}%
\def\Jtest{ }\def\Jstr{Duke Math. J.}%
\def\Vtest{ }\def\Vstr{66}%
\def\Dtest{ }\def\Dstr{1992}%
\def\Ptest{ }\def\Pstr{1--41}%
\def\Astr{\Underlinemark}%
\Refformat\egroup%

\bgroup\Resetstrings%
\def\Ecnt{0}\def\acnt{0}%
\def\Ftest{ }\def\Fstr{30}%
\def\Atest{ }\def\Astr{A\Initper \Initgap J\Initper  Silberger}%
\def\Ttest{ }\def\Tstr{The Knapp-Stein dimension theorem for $p$\snug-adic
groups}%
\def\Jtest{ }\def\Jstr{Proc. Amer. Math. Soc.}%
\def\Vtest{ }\def\Vstr{68}%
\def\Dtest{ }\def\Dstr{1978}%
\def\Ptest{ }\def\Pstr{243--246}%
\Refformat\egroup%

\bgroup\Resetstrings%
\def\Ecnt{0}\def\acnt{0}%
\def\Ftest{ }\def\Fstr{31}%
\def\Atest{ }\def\Astr{A\Initper \Initgap J\Initper  Silberger}%
\def\Ttest{ }\def\Tstr{Introduction to Harmonic Analysis on Reductive
$p$\snug-adic Groups}%
\def\Itest{ }\def\Istr{Princeton University Press}%
\def\Ctest{ }\def\Cstr{Princeton, NJ}%
\def\Stest{ }\def\Sstr{Mathematical Notes}%
\def\Ntest{ }\def\Nstr{23}%
\def\Dtest{ }\def\Dstr{1979}%
\def\Astr{\Underlinemark}%
\Refformat\egroup%

\bgroup\Resetstrings%
\def\Ecnt{0}\def\acnt{0}%
\def\Ftest{ }\def\Fstr{32}%
\def\Atest{ }\def\Astr{A\Initper \Initgap J\Initper  Silberger}%
\def\Ttest{ }\def\Tstr{The Knapp-Stein dimension theorem for $p$\snug-adic
groups. Correction}%
\def\Jtest{ }\def\Jstr{Proc. Amer. Math. Soc.}%
\def\Vtest{ }\def\Vstr{76}%
\def\Dtest{ }\def\Dstr{1979}%
\def\Ptest{ }\def\Pstr{169--170}%
\def\Astr{\Underlinemark}%
\Refformat\egroup%

\bgroup\Resetstrings%
\def\Ecnt{0}\def\acnt{0}%
\def\Ftest{ }\def\Fstr{33}%
\def\Atest{ }\def\Astr{A\Initper \Initgap V\Initper  Zelevinsky}%
\def\Ttest{ }\def\Tstr{Induced representations of reductive $p$\snug-adic
groups II, on irreducible representations of $GL(n)$}%
\def\Jtest{ }\def\Jstr{Ann. Sci. \'Ecole Norm. Sup. (4)}%
\def\Vtest{ }\def\Vstr{13}%
\def\Dtest{ }\def\Dstr{1980}%
\def\Ptest{ }\def\Pstr{165--210}%
\Refformat\egroup%

\endRefs
\enddocument
\end